\documentclass[11pt, reqno]{amsart}
\usepackage{mathrsfs}
\usepackage[active]{srcltx}
\usepackage{mathrsfs}
\usepackage{longtable}

\usepackage{xcolor}
\definecolor{bluecite}{HTML}{0875b7}
\usepackage[unicode=true,
bookmarksopen={true},
pdffitwindow=true,
colorlinks=true,
linkcolor=bluecite,
citecolor=bluecite,
urlcolor=bluecite,
hyperfootnotes=false,
pdfstartview={FitH},
pdfpagemode= UseNone]{hyperref}

\newcommand{\ds}{\displaystyle}

\newtheorem{example}{ Example}[section]
\newtheorem{proposition}{Proposition}[section]
\newtheorem{theorem}{Theorem}[section]
\newtheorem{lemma}{Lemma}[section]

\newtheorem{remark}{Remark}[section]
\numberwithin{equation}{section}

\usepackage{bbm}  
\usepackage{dsfont}

\begin{document}

\title[Sharp isoperimetric inequalities on manifolds]
{Sharp isoperimetric and Sobolev inequalities in spaces with nonnegative Ricci curvature
}

\author{Zolt\'an M. Balogh and Alexandru Krist\'aly}	

\thanks{Z. M. Balogh was
	supported by the Swiss National Science Foundation, Grant Nr. {200020\_191978}.  A. Krist\'aly  was supported by the UEFISCDI/CNCS grant PN-III-P4-ID-PCE2020-1001.}

\address{Mathematisches Institute,
	Universit\"at Bern,
	Sidlerstrasse 5,
	3012 Bern, Switzerland}

\email{zoltan.balogh@math.unibe.ch}

\address{Department of Economics, Babe\c s-Bolyai University, Str. Teodor Mihali 58-60, 400591 Cluj-Napoca,
	Romania \& Institute of Applied Mathematics, \'Obuda
	University, B\'ecsi \'ut 96/B, 1034
	Budapest, Hungary}

\email{alexandru.kristaly@ubbcluj.ro; kristaly.alexandru@nik.uni-obuda.hu.}



\begin{abstract}
	{\footnotesize \noindent 
		By using optimal mass transport theory we prove a sharp isoperimetric inequality in ${\sf CD} (0,N)$  metric measure spaces  assuming an asymptotic volume growth  at infi\-nity.  Our result extends recently proven  isoperimetric inequalities for normed spaces and Riemannian manifolds to a nonsmooth framework.   In the case of $n$-dimensional Riemannian manifolds with nonnegative Ricci curvature,  we outline an alternative proof of the rigidity result of S. Brendle (2021).
		 As applications of the isoperimetric inequality, we establish  Sobolev and Rayleigh-Faber-Krahn inequalities with explicit sharp constants  in Riemannian manifolds with nonnegative Ricci curvature; here we use appropriate symmetrization techniques and optimal volume non-collapsing properties.  
		The equality cases in the latter inequalities are also characte\-ri\-zed  by stating  that sufficiently smooth, nonzero extremal functions exist if and only if the Riemannian manifold is isometric to the Euclidean space.
	}
\end{abstract}
\maketitle



\section{Introduction}

The classical method in proving isoperimetric and related Sobolev inequalities in the Euc\-lidean space is by symmetrization arguments, see e.g. Talenti \cite{Talenti, Talenti-2} and Lieb and Loss \cite{LL}. Recent developments in this area showed two new approaches. The first one is highlighted by the works of  Cabr\'e  and  Ros-Oton \cite{Cabre-Ros-Oton} and Cabr\'e, Ros-Oton and Serra  \cite{Cabre-Ros-Oton-Serra} which is based on PDE techniques, notably on Aleksandrov-Bakelman-Pucci estimates; this is the so called ABP-method. The second approach is based on the theory of  optimal mass transportation (OMT, for short) and has been initiated in the seminal paper by Cordero-Erausquin,  Nazaret and Villani \cite{CE-N-Villani} and the monograph of Villani \cite{Villani}. 



Since the OMT-theory has been worked out in a rather general setting of  metric measure spaces due to Lott and Villani \cite{LV} and  Sturm \cite{Sturm-1, Sturm-2}, it is a natural question to ask  if this method is applicable to prove sharp geometric inequalities in general non-Euclidean settings. 

A pioneering result in this direction is due to Cordero-Erausquin, McCann and Schmuckenschl\"{a}ger \cite{CEMS} who extended the Borell-Brascamp-Lieb inequalities to Riemannian manifolds by using the OMT-theory; in fact, this work stood at the basis of the papers of \cite{LV, Sturm-1, Sturm-2}. 
 
In the case of Riemannian and weighted Riemannian manifolds satisfying a strictly positive lower bound on the Ricci curvature, sharp  L\'evy-Gromov-type isoperimetric inequalities have been obtained by Milman \cite{Milman, Milman2}.  An extension of Milman's results is due to Cavalletti and Mondino \cite{Cavalletti-Mondino}, who used  the OMT-approach in combination to the measure disintegration method of Klartag \cite{Klartag}  to prove a  sharp L\'evy-Gromov-type isoperimetric inequality on metric measure spaces satisfying the ${\sf CD}(\kappa,N)$ condition with $\kappa>0$ and $N\in  [1,\infty)$. Another extension of  \cite{Milman, Cavalletti-Mondino} is given by Ohta \cite{Ohta-1, Ohta-2} to the setting of (not necessarily reversible) Finsler manifolds.


Our goal is to show that the OMT-theory can also be successfully applied to prove sharp isoperimetric inequalities on  spaces with nonnegative Ricci curvature in a synthetic or classical sense. It is well known that without an additional condition no isoperimetric inequality holds true in general spaces with nonnegative Ricci curvature.  
 
For our results, the additional assumption is a volume growth control at infinity. To be more precise, let $N>1$ be a real number, and $(M,d,{\sf m})$ be a metric measure space satisfying the ${\sf CD}(0,N)$ condition, see  \cite{LV, Sturm-1, Sturm-2}. Let $B_x(r)=\{y\in M:d(x,y)<r\}$ be the metric ball with center $x\in M$ and radius $r>0.$ By the generalized Bishop-Gromov volume growth inequality, see Sturm \cite[Theorem 2.3]{Sturm-2}, it follows that 
$r\mapsto \frac{{\sf m}(B_x(r))}{r^N}$ is nonincreasing on $(0,\infty)$ for every $x\in M$. Moreover,  the \textit{asymptotic volume ratio} 
$${\sf AVR}_{M,d,{\sf m}}=\lim_{r\to \infty}\frac{{\sf m}(B_x(r))}{\omega_Nr^N}$$
is independent of the choice of $x\in M$.  (Here, the constant $\omega_N=\pi^{N/2}/\Gamma(1+N/2)$ plays  a normalization role and it is the volume of the Euclidean unit ball in $\mathbb R^N$ whenever $N\in \mathbb N$.)
Our standing assumption  is 
$$
{\sf AVR}_{M,d,{\sf m}}>0,
$$ and we say that  $(M,d,{\sf m})$ has \textit{Euclidean volume growth}.
In order to state our first result, we recall that the \textit{Minkowski content} of $\Omega\subset M$ is given by
$${\sf m}^+(\Omega)=\liminf_{\varepsilon \to 0^+}\frac{{\sf m}(\Omega_\varepsilon\setminus \Omega)}{\varepsilon},$$
where $\Omega_\varepsilon=\{x\in M:\exists y\in \Omega\ {\rm such\ that}\ d(x,y)< \varepsilon\}$ is the $\varepsilon$-neighborhood of $\Omega$ with respect to the metric $d$. 

\begin{theorem}\label{main-isop-CD} Let $(M,d,{\sf m})$ be a metric measure space satisfying the ${\sf CD}(0,N)$ condition for some $N>1,$ and having Euclidean volume growth. Then for every bounded Borel subset $\Omega\subset M$ it holds  
	\begin{equation}\label{eqn-isoperimetric}
		{\sf m}^+(\Omega)\geq N\omega_N^\frac{1}{N}{\sf AVR}_{M,d,{\sf m}}^\frac{1}{N}{\sf m}(\Omega)^\frac{N-1}{N}.	
	\end{equation}
	Moreover, inequality \eqref{eqn-isoperimetric} is sharp.
\end{theorem}

 Relying on the assumption of the Euclidean volume growth, our proof of \eqref{eqn-isoperimetric} uses careful limiting procedure in an OMT-based distorted Brunn-Minkowski inequality given by  Sturm \cite{Sturm-2}. Concerning the \textit{equality} in \eqref{eqn-isoperimetric},  we could expect certain rigidity statement in the spirit of Cavalletti and Mondino  \cite{Cavalletti-Mondino} who proved that if equality occurs in the isoperimetric inequality for some space $M$ with a strictly positive lower bound on its Ricci curvature, then $M$ is  isometric to a spherical suspension.
 At present time no such characterization is available for the equality in  \eqref{eqn-isoperimetric} in generic ${\sf CD}(0,N)$ spaces; the reason is that the argument based on the distorted Brunn-Minkowski inequality  - in spite of the fact that it provides the sharp inequality \eqref{eqn-isoperimetric} - seems to be too robust to identify the equality cases. 
 
Let us note that Theorem \ref{main-isop-CD} is in fact a generalization of the sharp isoperimetric inequalities  on \textit{weighted normed cones/spaces}, see Cabr\'e,  Ros-Oton and Serra \cite{Cabre-Ros-Oton-Serra}, on \textit{Riemannian manifolds}, see Brendle \cite{Brendle} and \textit{weighted Riemannian manifolds}, see  Johne \cite{Johne}.

Considering   a noncompact, complete  $n$-dimensional  Riemannian manifold $(M,g)$ with nonnegative Ricci curvature (${\sf Ric}\geq 0$, for short), this can be seen as a metric measure space $(M,d_g,{\rm Vol}_g)$, where $d_g$ and ${\rm Vol}_g$ denote the natural metric and canonical measure  on $(M,g),$ respectively.
 The  asymptotic volume
 ratio  of $(M,g)$ is a global geometric invariant  given by 
 ${\sf AVR}_g:={\sf AVR}_{M,d_g,{\rm Vol}_g}.$ By the Bishop-Gromov volume comparison principle one has that ${\sf
 	AVR}_g\leq 1$;  moreover,  ${\sf
 	AVR}_g= 1$ if and only if $(M,g)$ is isometric to the usual Euclidean space $(\mathbb R^n,g_0).$ Under the assumption of  Euclidean volume growth, i.e., $0< {\sf AVR}_g \leq 1$, Theorem \ref{main-isop-CD} implies that for every bounded and open subset $\Omega\subset M$ with $\mathcal{C}^1$ smooth
 boundary $\partial \Omega$, one has
 \begin{equation}\label{eq-isoperimetric-1}
 	\mathcal P_g(\partial \Omega)\geq  n\omega_n^\frac{1}{n} \ {\sf AVR}_g^\frac{1}{n}{\rm Vol}_g(\Omega)^\frac{n-1}{n},	
 \end{equation}
 where $\mathcal P_g(\partial \Omega)$ stands for the perimeter of $\partial \Omega$.
 
 We notice that \eqref{eq-isoperimetric-1} has been recently obtained by Brendle \cite{Brendle} using the ABP-method, who also proved that equality holds in \eqref{eq-isoperimetric-1}  for some $\Omega\subset M$ with $\mathcal{C}^1$ regular boundary if and only if ${\sf AVR}_g = 1$ and $\Omega$ is isometric to a ball $B\subset \mathbb R^n$. In our paper we also sketch an alternative proof of the characterization of the equality case based on the OMT-method.

Having inequality \eqref{eq-isoperimetric-1} at  hand - which is equivalent to a sharp $L^1$-Sobolev inequality on Riemannian manifolds with ${\sf Ric}\geq 0$ -, together with the characterization of the equality case, it is a natural question to consider the  validity of \textit{sharp $L^p$-Sobolev inequalities} in the same geometric setting whenever $p>1$. This problem has its genesis in the work of Aubin \cite{Aubin} who initiated in the early seventies  the determination of sharp constants in Sobolev inequalities in curved settings.
We notice that Aubin's program is rather well-understood in the counterpart setting of negative curvature. Indeed, sharp $L^p$-Sobolev inequalities hold  on Hadamard manifolds (i.e., simply connected, complete Riemannian manifold with nonpositive sectional curvature), see e.g.\ Muratori and Soave \cite{MS}, whenever the Cartan-Hadamard conjecture holds,  the latter being precisely the sharp Euclidean-type isoperimetric inequality on $(M,g)$. The validity of this conjecture is confirmed in low dimensions; see Weil \cite{Weil} in 2-dimension, Kleiner \cite{Kleiner} in 3-dimension and Croke \cite{Croke} in 4-dimension, respectively. 


For an  $n$-dimensional complete
Riemannian manifold $(M,g)$ with $n\geq 2$ and ${\sf Ric}\geq 0$,  endowed by its canonical measure ${\rm Vol}_g$, the simplest $L^p$-Sobolev inequality on $(M,g)$ reads as
$$
{ \left(\int_{M} |u|^{p^\star}{\rm d}v_g\right)^{1/p^\star}\leq
	C
	\left(\int_{M} |\nabla_g u|^p {\rm d}v_g\right)^{1/p}},\ \forall
u\in \mathcal C_0^\infty(M),\eqno{(\textbf{S})}
$$
where $C=C(n,p)>0$ is a universal constant, $1<p<n$ and $p^\star=\frac{pn}{n-p}$ is the critical Sobolev exponent. A local chart analysis shows that the validity of (\textbf{S}) necessarily implies  $C\geq {\sf AT}(n,p)$, where 
$${\sf AT}(n,p)=\pi^{-\frac{1}{2}}n^{-\frac{1}{p}} \left(\frac{p-1}{n-p}\right)^{1-1/p}\left(\frac{\Gamma(1+n/2)\Gamma(n)}{\Gamma(n/p)\Gamma(1+n-n/p)}
\right)^{{1}/{n}}$$
is the best Sobolev constant in $({\bf S})$ for the Euclidean space $(\mathbb R^n,g_0)$, see Aubin \cite{Aubin} and Talenti \cite{Talenti}. Moreover, in the Euclidean case,  the unique family of extremals is also identified.

Let us note first  that the existing literature already contains \textit{rigidity} results concerning inequality $({\bf S})$ on Riemannian manifolds  with ${\sf Ric}\geq 0$. Indeed, 
Ledoux \cite{Ledoux-CAG} proved that if $({\bf S})$ holds with $C={\sf AT}(n,p)$, then $(M,g)$ is isometric to the  Euclidean space $(\mathbb R^n,g_0).$ Moreover, a quantitative form of Ledoux's result were given  by do Carmo and Xia \cite{doCarmo-Xia} who proved that $(M,g)$ is topologically close to $(\mathbb R^n,g_0)$ whenever the constant $C>{\sf AT}(n,p)$ in $({\bf S})$ is sufficiently close to ${\sf AT}(n,p)$. 
In fact, a byproduct of do Carmo and Xia's approach is that the validity of $({\bf S})$ with a constant $C>0$ implies the  volume non-collapsing property 
\begin{equation}\label{volumenoncoll}
	{\sf AVR}_g\geq \left(\frac{{\sf AT}(n,p)}{C}\right)^n.
\end{equation}


\noindent  In particular, inequality \eqref{volumenoncoll} implies that for any complete Riemannian manifold $(M,g)$ with ${\sf Ric}\geq 0$ supporting the Sobolev inequality $({\bf S})$, one necessarily has that ${\sf
	AVR}_g>0,$ i.e., $(M,g)$ has Euclidean volume growth. 
The converse is also true that follows by a general result of  Coulhon and Saloff-Coste \cite{Coulhon-SC} where the constant $C>0$ in $(\bf S)$ is generically determined. 



%
%
%

Keeping our geometric setting, i.e.,\ $(M,g)$ is a noncompact, complete Riemannian manifold with ${\sf Ric}\geq 0$, our second  purpose is to provide \textit{sharp Sobolev inequalities on} $(M,g)$ by the sharp isoperimetric inequality from relation \eqref{eq-isoperimetric-1}. 
As expected, the asymptotic volume
ratio ${\sf
	AVR}_g\in (0,1]$ is explicitly encapsulated in these Sobolev inequalities and more spectacularly, they do provide sharp Sobolev constants. 

To formulate this result let us denote by   $\dot W^{1,p}(M)=\{u\in L^{p^\star}(M): |\nabla_g u|\in L^p(
M)\}.$  Using this notation we can state a sharp \textit{$L^p$-Sobolev inequality} in the spirit of ({\bf S}):



\begin{theorem}\label{main-0} Let $(M,g)$ be a noncompact, complete  $n$-dimensional Riemannian manifold $(n\geq 2)$ with ${\sf Ric}\geq 0$ having Euclidean volume growth, i.e., $0< {\sf AVR}_g \leq 1$.  If $p\in (1,n)$, then  for every $ u\in \dot W^{1,p}(M)$ one has
	\begin{equation}\label{egyenlet-0}
		\displaystyle  { \left(\int_{M} |u|^{p^\star}{\rm d}v_g\right)^{1/p^\star}\leq
			{\sf S}_g
			\left(\int_{M} |\nabla_g u|^p {\rm d}v_g\right)^{1/p}},
	\end{equation}
	where the constant $ {\sf S}_g={\sf AT}(n,p)\, {\sf AVR}_g^{-\frac{1}{n}}$ is   sharp. 
	Moreover,  equality holds in \eqref{egyenlet-0} for some nonzero and nonnegative function $u\in \mathcal{C}^{n}(M)\cap \dot W^{1,p}(M)$
	if and only if ${\sf AVR}_g=1,$ and thus $(M, g)$ is isometric to the Euclidean space $(\mathbb R^n, g_0)$.
\end{theorem}

As a byproduct, Theorem \ref{main-0} immediately implies the result of do Carmo and Xia \cite{doCarmo-Xia}. Indeed, if $({\bf S})$ holds for some $ C>0$, then by the sharpness of ${\sf S}_g$ we have $ C\geq {\sf S}_g,$ which is equivalent to \eqref{volumenoncoll}. Moreover, if $C={\sf AT}(n,p)$ then we obtain ${\sf AVR}_g\geq 1$, thus ${\sf AVR}_g= 1$, recovering Ledoux's rigidity result \cite{Ledoux-CAG} as well.

Note that inequality \eqref{egyenlet-0} belongs to the larger class of $L^p$-Gagliardo-Nirenberg inequalities $(1<p<n)$; we are going to treat such an inequality in Theorem \ref{main-GN} whose proof indicates the way to obtain further sharp Sobolev inequalities (i.e., $L^p$-log-Sobolev and $L^p$-Faber-Krahn  inequalities). 
In order to prove Theorem \ref{main-0}, we combine the sharp isoperimetric inequality  \eqref{eq-isoperimetric-1} with a {\it symmetrization argument from $(M,g)$ to $(\mathbb R^n,g_0)$}  in the spirit of Aubin \cite{Aubin}. In this way we establish a P\'olya-Szeg\H o-type inequality involving the number ${\sf AVR}_g$, see Proposition \ref{PSz-proposition}. Another challenge is to prove the \textit{sharpness} of the aforementioned inequalities (see e.g. the constant ${\sf S}_g$ in \eqref{egyenlet-0}). It turns out that  the  optimal volume non-collapsing properties of $(M,g)$  established by Ledoux \cite{Ledoux-CAG} (and by one of us \cite{Kristaly-Calculus}) provide precisely the required tool.

Another class of problems concerns the \textit{Rayleigh-Faber-Krahn} inequality where the sharp isoperimetrc inequality \eqref{eq-isoperimetric-1} provides  again the  strongest possible rigidity  statement. To formulate this result,  recall that given a complete $n$-dimensional Riemannian manifold $(M,g)$ (with no curvature restriction for the moment), it is well known that the first eigenvalue of the Beltrami-Laplace operator $-\Delta_g$ for the Dirichlet problem on a smooth bounded open set $\Omega\subset M$ has the  variational  characterization 
\begin{equation}\label{faber-krahn-0}
	\ds\lambda_{1,g}^D(\Omega)=\inf_{u\in \mathcal C_0^\infty(\Omega)\setminus \{0\}}\frac{\ds\int_{\Omega}|\nabla_g u|^2{\rm d}v_g}{\ds\int_{\Omega} u^2{\rm d}v_g}.
\end{equation}
According to Carron \cite{Carron2} (see also Hebey \cite[Proposition 8.1]{Hebey}), if $n\geq 3$ and ${\rm Vol}_g(M)=+\infty$, the validity of the general Sobolev inequality $(\bf S)$ is equivalent to the validity of a generic {Rayleigh-Faber-Krahn inequality} on $(M,g)$, i.e., there exists $\Lambda>0$ such that for any smooth bounded open set $\Omega\subset M$ one has
\begin{equation}\label{RFK-0}
	\lambda_{1,g}^D(\Omega)\geq \Lambda\, {\rm Vol}_g(\Omega)^{-\frac{2}{n}}.	
\end{equation}
In particular,  Theorem \ref{main-0} implies that inequality \eqref{RFK-0} holds  whenever $(M,g)$ is a  Riemannian  manifold with  ${\sf Ric}\geq 0$ having Euclidean volume growth $0< {\sf AVR}_g \leq 1$. In fact, in the latter geometric setting we can establish the \textit{sharp} form of \eqref{RFK-0}  (hereafter, $j_\nu$ stands for the first positive root of the Bessel function $J_\nu$ of the first kind with degree $\nu\in \mathbb R$): 
%

\begin{theorem}\label{RFK-1-0} Let $(M,g)$ be an $n$-dimensional Riemannian manifold as in Theorem \ref{main-0}.  Then for 	every bounded and open subset $\Omega\subset M$ with smooth
	boundary, we have 
	\begin{equation}\label{egyenlet-poincare-1-0}
		\displaystyle  \lambda_{1,g}^D(\Omega)\geq  {{\sf \Lambda}_g}\, {\rm Vol}_g(\Omega)^{-\frac{2}{n}},
	\end{equation}
	where the constant 
	$
	{\sf \Lambda}_g= j^2_{\frac{n}{2}-1}({\omega_n\, \sf AVR}_g)^{\frac{2}{n}}
	$
	is    sharp.\  Furthermore, equality holds in \eqref{egyenlet-poincare-1-0} for some bounded open subset $\Omega\subset M$ with smooth boundary if and only if ${\sf AVR}_g=1$ and $\Omega$ is isometric to a ball $B\subset \mathbb R^n$. 
\end{theorem}

The proof of Theorem \ref{RFK-1-0} is based on inequality \eqref{eq-isoperimetric-1} and fine properties of Bessel functions. We note that, while a similar result is recently established by Fogagnolo and Mazzieri \cite[Theorem 5.5]{Fogagnolo-Mazzieri} for $n\in \{3,...,7\}$,  Theorem \ref{main-0} is valid in any dimension $n\geq 2$.

The paper is organized as follows. Section \ref{section-sharp-CD} is devoted to sharp isoperimetric inequalities on ${\sf CD}(0,N)$ spaces. First, in \S \ref{subsection-2-1} we provide the proof for Theorem \ref{main-isop-CD} and present some  examples of ${\sf CD}(0,N)$ where our result applies. In \S \ref{section-OMT}, by using OMT-arguments, we outline a short, alternative proof to Brendle's rigidity result concerning the equality in \eqref{eq-isoperimetric-1}  in the context of Riemannian manifolds.    
In Section \ref{section-GN} we establish a sharp  Gagliardo-Nirenberg inequality on Riemannian manifolds  with ${\sf Ric}\geq 0$ having Euclidean volume growth, see Theorem  \ref{main-GN},   whose particular case is precisely Theorem  \ref{main-0}. To do this, we first provide a symmetrization argument, by establishing via the sharp isoperimetric inequality \eqref{eq-isoperimetric-1} an ${\sf AVR}_g$-dependent P\'olya-Szeg\H o inequality. 
 Theorem \ref{RFK-1-0} is proved in \S \ref{Poincare-section} by using the ${\sf AVR}_g$-dependent P\'olya-Szeg\H o inequality and fine features of Bessel functions.



\section{Sharp isoperimetric inequalities in ${\sf CD}(0,N)$ spaces 
}\label{section-sharp-CD}

\subsection{Proof of the sharp isoperimetric inequality}\label{subsection-2-1}

In this subsection we are going to prove inequality \eqref{eqn-isoperimetric}, its sharpness and provide some relevant examples and consequences. We first briefly recall the synthetic notion of nonnegative Ricci curvature introduced by Lott and Villani \cite{LV} and Sturm \cite{Sturm-1, Sturm-2}. 

Let $({M},d,\textsf{m})$ be a metric measure space, i.e.,
$({M},d)$ is a complete separable metric space, 
$\textsf{m}$ is a locally finite measure on $M$ endowed with its
Borel $\sigma$-algebra, and geodesic (for every two points  $x,y\in M$ there exists a minimizing geodesic $\gamma:[0,1]\to M$ parametrized proportional to arclength and $\gamma(0)=x$ and $\gamma(1)=y$). In particular, we can define the $s$-\textit{interpolant set} $Z_s(\cdot,\cdot)$, i.e., for every $(x,y) \in M \times M$,
\begin{eqnarray}\label{Riemannian-intermediate}
	Z_s(x,y) = \{ z \in M : d(x,z) = s d(x,y),\
	d(z,y) = (1-s) d(x,y)\},
\end{eqnarray}
and  for any nonempty sets $A,B \subset M$,  $$
Z_s(A,B) = \bigcup\limits_{(x,y) \in A \times B} Z_s(x,y).
$$ 
We assume that the measure
$\textsf{m}$ on $M$ is strictly positive, i.e.,
supp$[\textsf{m}]=M.$ As usual, $P_2(M,d)$ is the
$L^2$-Wasserstein space of probability measures on $M$, while
$P_2(M,d,\textsf{m})$ will denote the subspace of
$\textsf{m}$-absolutely continuous measures.

For  $N\geq 1,$ the {\it R\'enyi entropy functional}
${\rm Ent}_N(\cdot|\textsf{m}):P_2(M,d)\to \mathbb R$ with
respect to the measure $\textsf{m}$ is defined by
\begin{equation}\label{entropy}
	{\rm Ent}_N(\mu|\textsf{m})=-\int_M \rho^{-\frac{1}{N}}{\rm d}\mu=-\int_M \rho^{1-\frac{1}{N}}{\rm d}{\sf m},
\end{equation}
the function $\rho$ being the density of $\mu^{\rm ac}$ in
$\mu=\mu^{\rm ac}+\mu^{\rm s}=\rho \textsf{m}+\mu^{\rm s}$, where $\mu^{\rm ac}$ and $\mu^{\rm s}$
represent the absolutely continuous and singular parts of $\mu\in
P_2(M,d),$ respectively.

Given $N\geq 1$, the \textit{curvature-dimension condition} ${\sf CD}(0,N)$
 states that for all $N'\geq N$ the functional ${\rm Ent}_{N'}(\cdot|{\sf m})$ is
convex on the $L^2$-Wasserstein space $P_2(M, d,{\sf m})$, i.e.,  for each
$\mu_0,\mu_1\in  P_2(M,{d},\textsf{m})$ there exists
 a geodesic
$\Gamma:[0,1]\to  P_2(M,{d},\textsf{m})$ joining
$\mu_0$ and $\mu_1$ such that for every $s\in [0,1]$, 
$${\rm Ent}_{N'}(\Gamma(s)|\textsf{m})\leq (1-s) {\rm Ent}_{N'}(\mu_0|\textsf{m})+s {\rm Ent}_{N'}(\mu_1|\textsf{m}).$$
An almost immediate consequence of the above definition is the (distorted) \textit{Brunn-Minkowski inequality} on the metric measure space $({M},d,\textsf{m})$ satisfying the
 ${\sf CD}(0,N)$ condition; in particular, if $A,B\subset M$ are two Borel sets such that ${\sf m}(A)\neq 0\neq  {\sf m}(B)$, then for every $s\in [0,1]$ and $N'\geq N$ one has
 \begin{equation}\label{Brunn-Minkowski-initial}
 		{\sf m}(Z_s(A,B)^\frac{1}{N'}\geq (1-s){\sf m}(A)^\frac{1}{N'}+s{\sf m}(B)^\frac{1}{N'},
 \end{equation}
see e.g. Sturm \cite[Proposition  2.1]{Sturm-2}.\\

{\it Proof of inequality \eqref{eqn-isoperimetric}.}
Let us recall our setting:  
$(M,d,{\sf m})$ is a metric measure space satisfying the ${\sf CD}(0,N)$ condition for some real number $N>1,$ having Euclidean volume growth, and let $\Omega\subset M$  be a bounded Borel set.

Let $x_0\in \Omega$ and $R>0$ be arbitrarily fixed and let     
 $d_0:={\rm diam}(\Omega)<\infty.$ If $s\in [0,1]$, we claim that  
\begin{equation}\label{Z-inclusion}
Z_s(\Omega,B_{x_0}(R))\subseteq \Omega_{s(d_0+R)},
	\end{equation}
where $\Omega_\varepsilon$ is the $\varepsilon$-neighborhood of $\Omega$ with $\varepsilon> 0$. 
Let $s>0$ (for $s=0$,  \eqref{Z-inclusion} is trivial).  Indeed, if $z\in Z_s(\Omega,B_{x_0}(R))$, then by definition, there exist $x\in \Omega$ and $y\in B_{x_0}(R)$ such that $d(z,x)=sd(x,y)$ and $ d(z,y)=(1-s)d(x,y)$. In particular, we have that $d(x,y)\leq d(x,x_0)+d(x_0,y)< d_0+R$, thus, 
${\rm dist}(z,\Omega)\leq d(z,x)=sd(x,y)< s(d_0+R),$
which proves the claim \eqref{Z-inclusion}. 

On the other hand, the Brunn-Minkowski inequality \eqref{Brunn-Minkowski-initial} implies for every $s\in [0,1]$ that 
$$
	{\sf m}(Z_s(\Omega,B_{x_0}(R)))^\frac{1}{N}\geq (1-s){\sf m}(\Omega)^\frac{1}{N}+s{\sf m}(B_{x_0}(R))^\frac{1}{N}.
$$
Then, by \eqref{Z-inclusion} and using the definition of the Minkowski content, we have 
\begin{eqnarray*}
	{\sf m}^+(\Omega)&=&\liminf_{\delta\to 0}\frac{{\sf m}(\Omega_\delta\setminus \Omega)}{\delta}=\liminf_{s\to 0}\frac{{\sf m}(\Omega_{s(d_0+R)})-{\sf m}(\Omega)}{s(d_0+R)}\\&\geq &\liminf_{s\to 0}\frac{{\sf m}(Z_s(\Omega,B_{x_0}(R)))-{\sf m}(\Omega)}{s(d_0+R)}\\&\geq &\liminf_{s\to 0}\frac{\left((1-s){\sf m}(\Omega)^\frac{1}{N}+s{\sf m}(B_{x_0}(R))^\frac{1}{N}\right)^N-{\sf m}(\Omega)}{s(d_0+R)}\\&=&N{\sf m}(\Omega)^\frac{N-1}{N}\frac{{\sf m}(B_{x_0}(R))^\frac{1}{N}-{\sf m}(\Omega)^\frac{1}{N}}{d_0+R}.
\end{eqnarray*}
Since the latter estimate is valid for every $R>0$, we can take the limit $R\to \infty$ at the right hand side, and by the definition of the asymptotic volume ratio we obtain that
$${\sf m}^+(\Omega)\geq N\omega_N^\frac{1}{N}{\sf AVR}_{M,d,{\sf m}}^\frac{1}{N}{\sf m}(\Omega)^\frac{N-1}{N},$$
which is precisely relation \eqref{eqn-isoperimetric}. \\

{\it Sharpness of \eqref{eqn-isoperimetric}.}
In order to show the sharpness of \eqref{eqn-isoperimetric}, we argue by contradiction. Assume that there is a  constant 
$$ C > N\omega_N^\frac{1}{N}{\sf AVR}_{M,d,{\sf m}}^\frac{1}{N}$$
such that for all bounded Borel subsets $\Omega \subset M$ with positive measure it holds
\begin{equation} \label{eqn-contrad}
{\sf m}^+(\Omega)\geq C {{\sf m}}(\Omega)^\frac{N-1}{N}. 
\end{equation} 
To obtain the desired contradiction, we choose $ \Omega= B_{x_0}(r) $, $r>0$, and observe  that $(B_{x_0}(r))_{\delta} \subseteq B_{x_0}(r +\delta) $ for every $\delta>0$. By the monotonicity of the function $r \mapsto \frac{{{\sf m}}(B_{x_0}(r))}{ r^N}$
one has that 
$$  \frac{{{\sf m}}(B_{x_0}(r+ \delta ))}{( r+ \delta)^N} \leq  \frac{{{\sf m}}(B_{x_0}(r))}{ r^N}.$$
Consequently, we obtain that 
\begin{eqnarray*}
	\frac{{{\sf m}}((B_{x_0}(r))_{\delta})- {{\sf m}}(B_{x_0}(r)) }{\delta} \leq \frac{{{\sf m}}(B_{x_0}(r+\delta))- {{\sf m}}(B_{x_0}(r))}{\delta}\leq{{\sf m}}(B_{x_0}(r)) \frac{1}{\delta} \left[\left(\frac{r+\delta }{r}\right)^N-1\right].	
\end{eqnarray*}
Letting $\delta \to 0$ in the above inequality,  it follows that
\begin{equation}\label{eqn-intermediate-isoper}
	{\sf m}^+( B_{x_0}(r)) \leq N  \left( \frac{{{\sf m}}(B_{x_0}(r))}{r^N}\right)^{\frac{1}{N}} {{\sf m}}(B_{x_0}(r))^{\frac{N-1}{N}}.	
\end{equation}
Using the definition of the asympotic volume ratio,  
by \eqref{eqn-contrad} and \eqref{eqn-intermediate-isoper} one has 
	$$C\leq \lim_{r \to \infty} \frac{{\sf m}^+( B_{x_0}(r))}{{{\sf m}}(B_{x_0}(r))^{\frac{N-1}{N}}}\leq N \omega_N^\frac{1}{N}{\sf AVR}_{M,d,{\sf m}}^\frac{1}{N} < C,$$
a contradiction, proving the sharpness of \eqref{eqn-isoperimetric}. \hfill $\square$\\

\begin{remark}\rm
	Let us mention that the usage of Brunn-Minkowski inequalities is a powerful tool to prove isoperimetric and other related geometric inequalities; for recent results in this direction we refer to Kolesnikov and Milman \cite{Kolesnikov-Milman}, and Milman and Rotem \cite{Milman-Rotem}.    However, the double limiting process in the above proof of  \eqref{eqn-isoperimetric}  (by taking first  $s\to 0$ and then $R\to\infty$)   conceals those fine information that are crucial to characterize the equality cases.  Therefore, the use  of the  Brunn-Minkowski inequality \eqref{Brunn-Minkowski-initial} in general seems to be too rough to establish rigidity statements. 
\end{remark} 

In the sequel  we present some classes of  ${\sf CD}(0,N)$ spaces where our results can be applied.

\begin{example} \rm (\textbf{Weighted cones}) \label{ex-1}
	 Let $\Sigma\subseteq \mathbb R^n$ be an open convex cone with vertex at the origin, and $H:\mathbb R^n\to [0,\infty)$ be a gauge function (i.e., symmetric, convex and positively homogeneous of degree one). We endow the space $ \mathbb R^n$ with the induced metric $d_H(x,y) = H(x-y)$. Let $w$ be a continuous function in $\overline \Sigma$, positive in $E$, and positively homogeneous of degree $\alpha \geq 0$ such that $w^\frac{1}{\alpha}$ is concave in $\Sigma$ whenever $\alpha>0$. 
	In particular, the concavity of $w^\frac{1}{\alpha}$ is equivalent to the fact that the triplet $(\Sigma,d_H,w\mathcal L^n)$ is a ${\sf CD}(0,n+\alpha)$ space, see e.g. Villani \cite{Villani}, which also means that 	${\sf Ric}_w^\alpha\geq 0$ on $\Sigma$, see \eqref{Ric0} below.   
Moreover, by the homogeneity properties of $H$ and $w$, it is easy to check that 
	$$ {\sf AVR}_{\Sigma, d_H, w\mathcal L^n} = \frac{\displaystyle\int_{B_{d_H}(1) \cap \Sigma} w}{\omega_n} > 0.$$
	In particular, Theorem \ref{main-isop-CD} implies the main result of Cabr\'e,  Ros-Oton and Serra \cite[Theorem 1.3]{Cabre-Ros-Oton-Serra}, i.e.,
	for every open set $\Omega\subset \mathbb R^n$ with enough smooth boundary and $\displaystyle\int_{\Omega\cap \Sigma}w<\infty$ one has  
	\begin{equation}\label{Cabre-isop}
		\frac{P_{w,H}(\Omega;\Sigma)}{\displaystyle\left(\int_{\Omega\cap \Sigma}w\right)^\frac{n+\alpha-1}{n+\alpha}}\geq \frac{P_{w,H}(W;\Sigma)}{\displaystyle\left(\int_{W\cap \Sigma}w\right)^\frac{n+\alpha-1}{n+\alpha}},
	\end{equation}
	where $P_{w,H}$ denotes the anisotropic weighted perimeter w.r.t. $w$ and 
$H$, and $W$ is the Wulff set  associated with the gauge $H$, i.e., $ W= \{ x \in \mathbb R^n: x\cdot \nu \leq H(\nu), \ \nu \in \mathbb S^{n-1} \},$ see Wulff \cite{Wulff}.
	
	 We notice that the right hand side of the isoperimetric inequality \eqref{Cabre-isop} can be easily rewritten by means of ${\sf AVR}_{\Sigma,d_H,w\mathcal L^n}$, coming from the homogeneity of $w$ and $H$.   Particular forms of \eqref{Cabre-isop} can be found also in Cabr\'e and Ros-Oton \cite{Cabre-Ros-Oton}. 
	 Let us note that the images of the Wulff set $W$ under  dilations are isoperimetric sets. Very recently,  Cinti,  Glaudo,  Pratelli,  Ros-Oton and Serra \cite{Cinti-etal} proved that all isoperimetric sets in \eqref{Cabre-isop} are of this form (possibly up to some translations).
	
\end{example}

\begin{example}\rm (\textbf{Weighted Riemannian manifolds}) \label{ex-2}
	Let $(M,g)$ be a noncompact, complete,  $n$-dimensional Riemannian manifold and consider the Bakry-\'Emery Ricci curvature on the Riemannian metric measure space $(M,g,w {\rm d}v_g)$ given by 
	\begin{equation}\label{Ric0}
	{\sf Ric}_w^\alpha:={\sf Ric}-D^2(\log w)-\frac{1}{\alpha}D \log w \otimes D \log w,	
	\end{equation}
	where $\alpha> 0$, and $Dw$ and $D^2w$ denote the differential and Hessian of a differentiable function $w:M\to \mathbb (0,\infty)$. It is known that if ${\sf Ric}_w^\alpha\geq 0$ then $(M,g,w {\rm d}v_g)$ is a  ${\sf CD}(0,N)$ space with $N=n+\alpha.$ In particular, if $\Omega\subset M$ is an open set with smooth boundary, then ${\sf m}^+(\Omega)=\displaystyle \int_{\partial \Omega}w$ and ${\sf m}(\Omega)=\displaystyle\int_{ \Omega}w$, where ${\sf m}=w {\rm d}v_g$ is the weighted measure on $(M,g)$. Thus, assuming that
	$$\Lambda_\alpha=\lim_{r\to \infty}\frac{\displaystyle\int_{ B_x(r)}w}{r^{n+\alpha}}>0,$$
	then  Theorem \ref{main-isop-CD} applies and we immediately obtain the sharp isoperimetric inequality 
	\begin{equation}\label{Bakry-Emery}
		\displaystyle \int_{\partial \Omega}w\geq (n+\alpha) \Lambda_\alpha^\frac{1}{n+\alpha} \left(\displaystyle\int_{ \Omega}w\right)^\frac{n+\alpha-1}{n+\alpha}.	
	\end{equation}
	Inequality \eqref{Bakry-Emery} has been recently obtained by Johne \cite{Johne}, extending the ABP-method used by Brendle \cite{Brendle} in the unweighted case $w=1$. The equality case is characterized by Brendle \cite{Brendle} whenever $w=1$; this fact will be discussed in more details in \S \ref{section-OMT}.
\end{example}

\begin{example}\rm \label{example-BS} (\textbf{Euclidean cones})
	Following \cite{BS} and  \cite{Ketterer},  we present here a  general class of examples such that we have equality in \eqref{eqn-isoperimetric}; these are the so-called {\it Euclidean cones} and are defined as follows. Let $(M, d, {\sf m})$ be a complete metric measure space such that ${\sf diam} M \leq \pi$. We define ${\sf Con} (M)$ as the quotient of $M\times [0, \infty)$ by identifying all points of $M\times\{0\}$ with a single point $O$ - the origin of ${\sf Con}(M)$. The metric on ${\sf Con}(M)$ is defined by 
	$$
		d_{c}((x,s), (y,t)) = \sqrt{s^2 + t^2 -2st\cos d(x,y)}, 
$$
	and the measure by $d{\sf m_c}(x,s) = d{\sf m}(x) \otimes s^n ds$. 
	Let us assume that the space $({\sf Con}(M), d_c, {\sf m_c})$ satisfies the ${\sf CD} (0, n+1) $ condition. 
	According to Bacher and Sturm \cite{BS}, this holds true when $M$ is an $n$-dimensional Riemannian manifold with ${\sf Ric}  \geq n-1$, or more generally, if $M$ is a weighted Riemannian manifold satisfying the 
	${\sf CD}(n-1, n)$ condition (even when $n$ is not necessarily an integer). Moreover, according to Ketterer \cite{Ketterer} the same holds if we consider $M$ to be a ${\sf CD^{\ast}}(n-1, n)$ space. 
A direct calculation gives that 
	$$ {\sf AVR}_{{\sf Con}(M), d_c, {\sf m_c}}  = \frac{{\sf m}(M)}{(n+1)\omega_{n+1}} >0 .$$
	In conclusion the statement of Theorem  \ref{main-isop-CD} applies. 
	
	We claim that balls $B_O(R)$ centered at the origin of the cone are isoperimetric sets for all $R>0$. To see this, first note that 
	$ {\sf m_c}(B_O(R))= \frac{{\sf m}(M)R^{n+1}}{n+1}, $ and furthermore we  claim that  $B_O(R+\epsilon) = (B_{O}(R))_{\epsilon}$. To see this last equality we note first that the 
	inclusion  $B_O(R+\epsilon) \subseteq  (B_{O}(R))_{\epsilon}$ is trivial. To check the opposite inclusion, pick a point $(y,t) \in (B_O(R))_{\epsilon}$; then there exists a point $(x,s) \in B_O(R)$ such that 
	$ \epsilon > d_{c}((x,s), (y,t)) = \sqrt{s^2 + t^2 -2st\cos d(x,y)} \geq |s-t|.$
	Since $(x,s) \in B_O(R)$ we have that $s < R$ and thus $t < R+ \epsilon$. From here we obtain that ${\sf m_c^+}(B_O(R)) = {\sf m}(M) R^n$ and we conclude that equality holds in \eqref{eqn-isoperimetric} for $B_O(R)$. We do not know if all isoperimetric sets are of this form.
	
\end{example}

\begin{remark}\rm The characterization of the equality case in \eqref{eqn-isoperimetric} is a challenging problem even in the particular settings of the above examples. In the case of  Example \ref{ex-1} 
(weighted cones), a careful stability argument is carried out in \cite{Cinti-etal} in order to 	characterize the 
 isoperimetric sets in \eqref{Cabre-isop}, while in the case of Example \ref{ex-2} we expect a strong rigidity of the manifold (and presumably of the weight) akin to the one stated by Brendle \cite{Brendle} (and discussed below) in the unweighted setting. 
The question of characterization of the equality case in the isoperimetric inequality  in Theorem  \ref{main-isop-CD} for general metric measure spaces requires further investigations that will be considered in a forthcoming work. 
\end{remark}

\subsection{Rigid isoperimetric inequalities in Riemannian manifolds with ${\sf Ric}\geq 0:$ the canonical case} 
\label{section-OMT}


In this subsection  we focus to  the following result: 

\begin{theorem}\label{tetel-egyenloseg} {\rm (Brendle \cite{Brendle})}   Let $(M,g)$ be an $n$-dimensional Riemannian manifold as in Theorem \ref{main-0}. Equality holds in   \eqref{eq-isoperimetric-1} for some $\Omega\subset M$ with $\mathcal{C}^1$ regular boundary if and only if ${\sf AVR}_g = 1$ and $\Omega$ is isometric to a ball $B\subset \mathbb R^n$.
\end{theorem}

Let us note first that besides Brendle's ABP-based proof (which is valid in any dimension), Theorem \ref{tetel-egyenloseg} has been proven in the 3-dimensional case by Agostiniani, Fogagnolo and Mazzieri \cite{AFM} by using Huisken's mean curvature
flows; moreover, Fogagnolo and Mazzieri \cite{Fogagnolo-Mazzieri} extended their arguments to manifolds up to 7 dimensions. We shall also outline a short, alternative proof of Theorem \ref{tetel-egyenloseg} by using tools from the OMT-theory. Our primordial motivation by doing so  is that we present an approach that might be useful also to wider classes of possibly nonsmooth settings. At this stage however,  certain technical issues prevent us to  carry out the proof in more general structures. \\
 

{\it Outline of the proof of Theorem \ref{tetel-egyenloseg}.} Let  $\Omega\subset M$ be a bounded, connected and open set with smooth
boundary such that equality holds in  \eqref{eq-isoperimetric-1}. 
We divide the proof in two steps.

{\bf Step 1:} \textit{we show that at the points of the isoperimetric set $\Omega$ the manifold $M$ is locally isometric to the Euclidean space}.
Let us consider the probability measures $$\mu=\frac{\mathds{1}_\Omega}{{{\rm Vol}_g}(\Omega)}{\rm d}v_g\ \ {\rm  and}\ \ \nu=\frac{\mathds{1}_{\Omega_r}}{{{\rm Vol}_g}(\Omega_r)}{\rm d}v_g$$ and the associated optimal transport map $T_r(x)=\exp_x(-\nabla_g u_r(x))$   for  a.e.\ $ x\in \overline \Omega,$ 
where   $\mathds{1}_{A}$ denotes the indicator function of the set $ A\subset M$, and $u_r:\overline \Omega\to \mathbb R$ is a $c=d_g^2/2$-concave function, see Cordero-Erausquin, McCann and  Schmuckenschl\"{a}ger \cite{CEMS} and McCann \cite{McCann}. Then $T_r: \overline \Omega \to \Omega_r$ is injective except of a null set and the change of variables formula holds; in particular, we have the Monge-Amp\`ere equation 
\begin{equation}\label{Monge-Ampere}
	\frac{1}{{{\rm Vol}_g}(\Omega)}=\frac{1}{{{\rm Vol}_g}(\Omega_r)}{\rm det}DT_r(x)\ \ {\rm for\ a.e.}\ x\in \Omega.	
\end{equation}
By the construction of the optimal transport map $T_r$, it turns out that $|\nabla_g u_r(x)|\leq r+2d_0$ for a.e. $x\in \overline \Omega$, where $d_0={\rm diam}(\Omega)$.  Let $t>0$ be arbitrarily fixed. For  $r\geq t$, we introduce the family of scaling functions $w_{t, r} = t  \frac{u_r}{r}.$ 
The latter estimate implies that there exists $C_0>0$ (not depending on $r$) such that for every $r\geq t$ and a.e. $x\in \overline \Omega$, 
$
	|\nabla_g w_{t,r}(x)| \leq tC_0,
$
i.e.,  the family $\{w_{t,r}\}_{r\geq t}$ is equicontinuous. 
 By Arzel\`a-Ascoli's theorem we obtain a sequence $\{w_{t, r_k}\}_k$ that converges uniformly to some function $w_t:\overline \Omega\to \mathbb R$ as $k \to \infty$ (and $r_k\to \infty$);  moreover,   $w_t$ is $c$-concave, which follows by the general theory of $c$-concave functions (see  Villani \cite{Villani}).
 
 Up to a smoothing argument \`a la Greene and Wu \cite{GW, GW1}, we assume that $w_{t,r}$ is enough regular, thus the divergence theorem, Schwarz inequality and the equality in  \eqref{eq-isoperimetric-1} yield  
 \begin{eqnarray}\label{first-estimate}
 	\nonumber	\int_\Omega \left(1-\frac{\Delta_g w_{t,r}(x)}{n}  \right) {\rm d}v_g&=&{\rm Vol}_g(\Omega)-\frac{1}{n}\int_{\partial \Omega} \langle \nabla_g w_{t,r},{\bf n} \rangle_g {\rm d}\sigma_g\\&\leq &\nonumber {\rm Vol}_g(\Omega)+\frac{1}{n}\frac{t}{r}(r+2d_0)\mathcal P_g(\partial \Omega)\\&=&{\rm Vol}_g(\Omega)+\frac{t}{r}(r+2d_0)\omega_n^\frac{1}{n}{\sf AVR}_g^\frac{1}{n}{{\rm Vol}_g}(\Omega)^\frac{n-1}{n},
 \end{eqnarray}
 where ${\bf n}(x)$ stands for the unit outward normal vector at $x \in \partial \Omega$.  Since ${\sf Ric}\geq 0$,  the volume distortions in
 $(M, g)$ verify
 \begin{equation}\label{vol-distortion}
 	v_s(x,y) = \lim\limits_{r \to 0}\frac{{\rm Vol}_g\left(Z_s(x, B_y( r))\right)}{{\rm Vol}_g\left(B_y(sr)\right)}\geq 1
 \end{equation}
 for every $s\in (0,1)$ and  $x,y\in M$ with $y\notin {\sf cut}(x)$,  
 see \cite[Corollary 2.2]{CEMS}, where ${\sf cut}(x)\subset M$ is the cut-locus of $x.$ Thus, the Jacobian determinant inequality from \cite[Lemma 6.1]{CEMS} becomes 
 \begin{eqnarray} \label{eq-interpolation2-0}
 \nonumber	\left( \det DT_{t,r} (x) \right)^{\frac{1}{n}} &\geq& \left(1-\frac{t}{r}\right)v_{1-t/r}(T_r(x),x)^\frac{1}{n} + \frac{t}{r}v_{t/r}(x,T_r(x))^{\frac{1}{n}} \left( \det DT_{r}(x)\right)^{\frac{1}{n}}\\&\geq& 1-\frac{t}{r} + \frac{t}{r} \left( \det DT_{r}(x)\right)^{\frac{1}{n}} \ {\rm for\ a.e.} \ x \in \Omega. 
 \end{eqnarray}
 Combining \eqref{Monge-Ampere}, \eqref{first-estimate} and \eqref{eq-interpolation2-0}, and letting $r\to \infty$, it follows by Fatou's lemma that  
 $$\int_\Omega \left(1-\frac{\Delta_g w_{t}}{n} - \left({\rm det}DF_{t} \right)^\frac{1}{n} \right) {\rm d}v_g  \leq 0,$$
 where $F_t: \overline\Omega  \to F_t(\overline\Omega)$ is the optimal transport map 
 $ F_t(x) = \exp_x (-\nabla_g  w_t(x)).$
 On the other hand, since $w_t:\overline \Omega\to \mathbb R$ is $c$-concave, one has the pointwise estimate
 \begin{equation} \label{eq:wang-zhang} 
 0\leq 1-\frac{\Delta_g w_{t}(x)}{n} - \left({\rm det}DF_{t} (x)\right)^\frac{1}{n}\ \ {\rm for\ a.e.}\ x\in {\Omega}, 
 \end{equation}
 see e.g. Wang and  Zhang \cite{WZ}. 
 By the latter two inequalities we obtain the second order PDE 
 
 \begin{equation} \label{eq-important}
 	1 - \frac{\Delta_g w_t(x)}{n}= \left( {\rm det} DF_t(x)\right)^{\frac{1}{n}} \ \ {\rm for\ a.e.}\ \ x\in \Omega.
 \end{equation}
Having  \eqref{eq-important} for every $t>0$, and by using that  $w_{st} = s w_t$, the Jacobian determinant inequality   for the map $F_{st}$ and relation  \eqref{vol-distortion} imply   for a.e. $ x\in \Omega$  that
\begin{eqnarray*}
	1- s\frac{\Delta_g w_{t}(x)}{n} &=&1- \frac{\Delta_g w_{st}(x)}{n} = \left({\rm det} DF_{st}(x)\right)^\frac{1}{n} \\
	&\geq & (1-s) \left( v_{1-s}(F_t(x), x)\right)^\frac{1}{n} + s \left(v_s(x, F_t(x))\right)^\frac{1}{n} \left({\rm det} DF_t(x)\right)^{\frac{1}{n}}
	\\ & \geq & 1-s + s\left(1-\frac{\Delta_g w_{t}(x)}{n} \right) = 1- s\frac{\Delta_g w_{t}(x)}{n}.
\end{eqnarray*}
According to these estimates, we must have equalities everywhere in the above chain of inequalities; thus, we necessarily have   that the volume  distortions should verify
\begin{equation} \label{eq-distort=1}
	v_{1-s}(F_t(x), x ) = v_s(x, F_t(x))= 1\ {\rm for\ a.e.}\  \ x\in \overline\Omega.
\end{equation}
By our earlier result \cite[Theorem 4.1]{BK-Adv} (see also  Chavel \cite[Theorem III.4.3]{Chavel}), we can conclude that the sectional curvatures  along the geodesic segments $s \mapsto \exp_x (-s\nabla_g  w_t(x))=F_{ts}(x)$, $s\in [0,1]$, connecting $x$ to $F_t(x)$ are constantly equal to $0$;  
this shows in particular that in all points of $\overline \Omega$ the sectional curvatures identically vanish. This fact implies that  $\overline \Omega \subset M$ is locally isometric to the Euclidean space, see e.g. Petersen \cite[Theorem 5.5.8]{Petersen}.  This concludes the first step of the proof. 

{\bf Step 2:} {\it we upgrade the local isometries to a global one.}  
In order to carry out the second step we use first a covering argument combined with the Bishop-Gromov comparison and \eqref{eq-important} to conclude that for some $\alpha>0,$
\begin{equation}\label{Hessian-identity}
	{\rm Hess}_g (-w_t)(x) = t\alpha {\rm Id}\ \ {\rm for\ every}\ t>0,\ x\in \Omega.
\end{equation}
In particular, we have that $\Delta_g(-w_t)=nt\alpha$, and by \eqref{eq-important}, it yields
\begin{equation}\label{diff-DFt}
	\det DF_t(x)=(1+\alpha t)^n\ \ {\rm for\ every}\ t>0,\ x\in \Omega.
\end{equation}

%

\noindent The Jacobian determinant inequality, the equality in \eqref{eq-isoperimetric-1} and  the Monge-Amp\`ere equation \eqref{Monge-Ampere} yield 
$$	\int_\Omega \left({\rm det}DF_{t} (x)\right)^\frac{1}{n} {\rm d}v_g\geq {\rm Vol}_g(\Omega)+\frac{t}{n}\mathcal P_g(\partial \Omega).$$
By \eqref{eq-important} and the homogeneity property $w_t=tw$,  the latter estimate becomes equivalent to 
$$\mathcal P_g(\partial \Omega)\leq \int_\Omega \Delta_g (-w) {\rm d} v_g.$$

\noindent 	Let us observe that $\Delta_g w = -n\alpha=$constant, thus $w$ is smooth up to the boundary $\partial \Omega$. Since $ |\nabla_g w | \leq 1$ on $\Omega$,  by the divergence theorem and Schwarz inequality  we obtain that
\begin{equation}\label{perimeter-uj}
	\mathcal{P}_g (\partial \Omega) \leq \int_\Omega \Delta_g (-w) {\rm d} v_g=\int_{\partial \Omega}\langle\nabla_g(-w)(x),  {\bf n}(x)\rangle_g{\rm d}\sigma_g(x)\leq \int_{\partial \Omega}{\rm d}\sigma_g=\mathcal{P}_g (\partial \Omega).
\end{equation}		
In particular, we have equality in the Schwarz inequality, which implies that 
\begin{equation}\label{boundary-value-condition}
	\nabla_g(-w)(x)={\bf n}(x)\ \ {\rm 
	for\ every}\ \ 
	x\in \partial \Omega.
\end{equation}

If $x_0\in \overline \Omega$ is the global minimum of $-w$ over $\overline \Omega$, by \eqref{boundary-value-condition} we clearly have that $x_0$ cannot belong to $\partial \Omega$; thus $x_0\in \Omega$ and subsequently, $\nabla_g w(x_0)=0$. Thus, $x_0\in  \Omega$ is a fixed point of $x\mapsto F_t(x)=\exp_x (-t\nabla_g  w(x))$ for every $t>0$. 

 We now show that the optimal transport map $F_t$ passes information from the infinitesimal volume distortion at the critical point $x_0$ of $w_t$ to infinity.  To do that, let $B_{x_0}(r)\subset \Omega$ be a ball with enough small radius $r>0$.  Since  $B_{x_0}(r)\subset \Omega$ is isometric to a ball in $\mathbb R^n$ with the same radius $r>0$, for every $t>0$ it follows by
\eqref{diff-DFt}  that  
\begin{equation}\label{Ft-vol}
	{{\rm Vol}_g}(F_t(B_{x_0}(r)))=\int_{B_{x_0}(r)}\det DF_t(x){\rm d}v_g=(1+\alpha t)^n {{\rm Vol}_g}(B_{x_0}(r))= (1+\alpha t)^n\omega_nr^n.
\end{equation}
Due to \eqref{Hessian-identity} and $\nabla_g w(x_0)=0$, a simple estimate shows that there exists $C>0$ such that  
for every $t>0$, one has 
$
	F_t(B_{x_0}(r))\subseteq B_{x_0}(t(\alpha r+Cr^2)+r).
$
Combining \eqref{Ft-vol} and the latter inclusion, we obtain 
$\omega_nr^n(1+\alpha t)^n \leq {{\rm Vol}_g}(B_{x_0}(t(\alpha r+Cr^2)+r)).$
For a fixed $r>0$, dividing by $t^n$ and taking  $t\to \infty$, the definition of ${\sf AVR}_g$ implies that 
$r^n\alpha^n \leq {\sf AVR}_g  (\alpha r+Cr^2)^n.$
Dividing by $r^n$ this inequality and taking $r\to 0$, it follows that $1\leq {\sf AVR}_g$. Thus, we have ${\sf AVR}_g=1$, i.e., $(M,g)$ is isometric to the Euclidean space $(\mathbb R^n,g_0)$. In particular, it follows that $\mathcal P_g(\partial \Omega)= n\omega_n^\frac{1}{n}{{\rm Vol}_g}(\Omega)^\frac{n-1}{n};$ 
being in the Euclidean setting (up to an isometry), the latter equality implies that $\Omega\subset M$ is isometric to a ball    $B\subset \mathbb R^n$. 
\hfill$\square$

\begin{remark}\rm 
	We notice the 'duality' of our OMT-argument with respect to Brendle's proof. On one hand, the ABP-method applied by Brendle \cite{Brendle} begins with a specific  PDE akin to \eqref{eq-important}  with a Neumann boundary value condition and uses an estimate of the type \eqref{eq:wang-zhang} to arrive via a Ricci flow to the rigidity result. On the other hand, our OMT-argument begins with optimal transport rays and the Monge-Amp\`ere equation \eqref{Monge-Ampere} by using \eqref{eq:wang-zhang} to conclude the PDE \eqref{eq-important} together with the Neumann boundary value condition \eqref{boundary-value-condition},  whose solution gives the required information about the whole manifold, i.e., ${\sf AVR}_g=1$.  
\end{remark}

\begin{remark}\rm 
	The smoothness of the boundary of the isoperimetric set is an essential requirement not only in our argument, but also in  
	Agostiniani,  Fogagnolo and Mazzieri \cite{AFM,Fogagnolo-Mazzieri} and Brendle \cite {Brendle}. However, we expect that the smoothness assumption on the boundary  might be replaced by a more general condition. Indeed, a careful  inspection of our proof shows that the same argument can be extended to cover the case of domains $\Omega$ with \textit{Lipschitz regular  boundaries}.
	On the other hand, isoperimetric sets $\Omega\subset M$  (i.e., satisfying equality in \eqref{eqn-isoperimetric} or \eqref{eq-isoperimetric-1}) are sets of finite perimeter. In the Euclidean case, the structure of sets with finite perimeter is well-understood, see   Ambrosio,  Fusco and Pallara   \cite{AFP}; this should give useful information on the regularity of  $\partial \Omega$ also in our case. Therefore, we believe that the rigidity statement should hold true with no additional apriori  boundary regularity assumption as such a property should already be encoded into the initial fact that $\Omega$ is a set of finite perimeter.  However, the proof of such a general statement is far from trivial, even in the Euclidean setting, see e.g.\ the survey notes by Fusco \cite{Fusco} on the early works of De Giorgi. 
\end{remark} 

Riemannian manifolds with ${\sf Ric}\geq 0$ have been widely studied in the literature, stating various classifications and topological rigidities, see e.g. Anderson \cite{Anderson}, Cheeger and  Colding \cite{CC}, Colding \cite{Colding}, Li \cite{Li}, Liu \cite{Liu}, Menguy \cite{Menguy}, Perelman \cite{Perelman-JAMS}, Reiris \cite{Reiris}, Zhu \cite{Zhu}. To conclude this section,  we present two Riemannian  manifolds with ${\sf Ric}\geq 0$ that satisfy in addition also the Euclidean volume growth condition,  providing as well their explicit asymptotic volume ratios. 

\begin{example}\rm (\textbf{Rotationally invariant metric on $\mathbb R^n$})  
Let $n\geq 3$ and   $f:[0,\infty)\to [0,1]$ be a smooth nonincreasing function such that
$f(0)=1$ and $\ds\lim_{s\to \infty} f(s)=a\in (0,1]$.  
We consider the \textit{rotationally invariant metric} on $\mathbb R^n$ defined by the warped product metric 
$$g={\rm d}r^2+F(r)^2 {\rm d}\theta^2,$$
where  $F(r)=\displaystyle\int_0^r f(s){\rm d}s$ and 
${\rm d}\theta^2$ is the standard metric on the sphere $\mathbb S^{n-1}$. If $x=(x_1,\theta_1)$ and $\tilde x=(x_2,\theta_2)$ are two points in $\mathbb R^n$, it turns out that $d_g(x,\tilde x)\geq |x_1-x_2|$, which implies that $(M,g)$ is complete. Furthermore, it is well known that 
the sectional (thus, the Ricci) curvature of $(\mathbb R^n,g)$ is nonnegative, 
see Carron \cite{Carron}.

For $R\gg 1,$ one has that
${\rm Vol}_g(B_0(R))=\displaystyle\int_{B_0(R)}{\rm d}v_g\sim n\omega_n\displaystyle\int_0^R F(r)^{n-1}{\rm d}r.$
The latter estimate and L'H\^ospital's rule give that 
$${\sf AVR}_g=\lim_{R\to \infty}\frac{{\rm Vol}_g(B_0(R))}{\omega_n R^n}=\lim_{R\to \infty}\frac{n\ds\int_0^R F(r)^{n-1}{d}r}{R^n}=\lim_{R\to \infty}\frac{F(R)^{n-1}}{R^{n-1}}=a^{n-1}\in (0,1].$$

When $a=1$, i.e., ${\sf AVR}_g=1$, by our monotonicity assumption it turns out that $f\equiv1$ on $[0,\infty)$; thus $F(r)=r$  and the metric    $g=g_0={\rm d}r^2+r^2 {\rm d}\theta^2$ becomes Euclidean.  



\end{example}

\begin{example}\rm ({\bf Asymptotically locally Euclidean manifolds}) Following  Agostiniani,  Fogagnolo and Mazzieri \cite[Definition 4.13]{AFM}, a complete, noncompact Riemannian manifold $(M,g)$ is {\it asymptotically locally Euclidean manifold} if  there exist a compact set $K\subset M$, a ball $B\subset \mathbb R^n$, a diffeomorphism $\Psi:M\subset K\to \mathbb R^n\setminus B$, a number $\tau>0$ and a finite subgroup $G$ of $SO(n)$ acting freely on $\mathbb R^n\setminus B$ such that
\begin{equation}\label{ALE-1}
	(\Psi^{-1}\circ \pi)^*g(z)=g_0+O(|z|)^{-\tau};
\end{equation} 
\begin{equation}\label{ALE-2}
	\left| \partial_i((\Psi^{-1}\circ \pi)^*g)\right|(z)=O(|z|)^{-\tau-1};
\end{equation} 
\begin{equation}\label{ALE-3}
	\left| \partial_i\partial_j((\Psi^{-1}\circ \pi)^*g)\right|(z)=O(|z|)^{-\tau-2},
\end{equation} 
where $\pi:\mathbb R^n\to \mathbb R^n/G$ stands for the natural projection, $z\in \mathbb R^n\setminus B$ and $i,j\in \{1,...,n\}.$

Due to assumptions \eqref{ALE-1}-\eqref{ALE-3}, it turns out that $(M,g)$ has Euclidean volume growth; furthermore, one has that
\begin{equation}\label{finite}
	{\sf AVR}_g=\frac{1}{{\rm Card}(G)},	
\end{equation}
see \cite[rel. (4.31)]{AFM}. 
In particular, $(M,g)$ is isometric to $(\mathbb R^n,g_0)$ if and only if $G=\{{\rm Id}\}\subset  SO(n);$ otherwise, $0<{\sf AVR}_g<1$. 

When $n=3$, the finite subgroups of $SO(3)$ are isomorphic to either  a  cyclic  group $\mathbb Z/m=\mathbb Z_m$ $(m\in \mathbb N\setminus \{0,1\})$,  a  dihedral  group $D_m$,  or  the  rotational symmetry group of a regular solid, i.e., (a) the symmetry group of the tetrahedron $A_4$, (b) the symmetry group of the cube $S_4$ (or octahedron), (c) the symmetry group of the dodecahedron $A_5$ (or icosahedron). These subgroups of $SO(3)$ together with \eqref{finite} can be efficiently applied to produce sharp isoperimetric inequalities on 3-dimensional asymptotically locally Euclidean manifolds. 
\end{example}

 \section{Sharp and rigid Sobolev inequalities on Riemannian manifolds with ${\sf Ric}\geq 0$} \label{section-GN}

%

 
 Let $u:M\to \mathbb R$ be a fast decaying function, i.e., ${\rm Vol}_g(\{x\in M:|u(x)|>t\})<+\infty$ for every $t>0.$ For such a function, let 
 inspired by  Aubin \cite{Aubin} and Druet, Hebey and Vaugon \cite{Druetetal}, we associate  its
 Euclidean rearrangement function $u^\star:\mathbb R^n\to
 [0,\infty)$ which is
 radially symmetric, nonincreasing in $|x|$, and for every $t>0$ is defined by
 \begin{equation}\label{vol-egyenloseg}
 	{\rm Vol}_{g_0}\left(\{x\in \mathbb R^n:u^\star(x)>t\}\right)={\rm Vol}_g\left(\{x\in
 	M:|u(x)|>t\}\right).
 \end{equation}
 By  \eqref{vol-egyenloseg} and the layer cake representation, see Lieb and Loss \cite{LL}, it turns out that 
 \begin{equation}\label{elso-tul}
 	{\rm Vol}_g({\rm supp}(u))={\rm Vol}_{g_0}({\rm
 		supp}(u^\star)),
 \end{equation} 
 and if $u\in L^q(M)$ for some $q\in (0,\infty)$, the Cavalieri principle reads as 
 \begin{equation}\label{masodik-tul}
 	\|u\|_{L^q(M)}=\|u^\star\|_{L^q(\mathbb R^n)}.
 \end{equation}

 
 For every $0<t<\|u\|_{L^\infty(M)}$, let us consider the sets 
 \begin{equation}\label{omegak}
 	\Omega_t=\{x\in M:|u(x)|>t\}
 	\ \ {\rm and}\ \ \Omega_t^\star=\{x\in \mathbb R^n:u^\star(x)>t\},
 \end{equation} 
 respectively. The key ingredient in our arguments is the following P\'olya-Szeg\H o inequality.

 
 \begin{proposition}\label{PSz-proposition} Let $(M,g)$ be an $n$-dimensional  Riemannian manifold with ${\sf Ric}\geq 0$ and having Euclidean volume growth, 
 	  and $u:M\to \mathbb R$ be a fast decaying function such that $|\nabla_g u|\in L^p(M)$, $p>1$. Then  one has 
 	\begin{equation}\label{harmadik-tul}
 		\|\nabla_g u\|_{L^p(M)}\geq {\sf AVR}_g^\frac{1}{n}\ \|\nabla u^\star\|_{L^p(\mathbb
 			R^n)}.
 	\end{equation} 
 	In addition, if equality holds in \eqref{harmadik-tul} for some nonnegative $u\in \mathcal{C}^{n}(M)\setminus \{0\}$, then ${\sf AVR}_g=1,$ i.e., $(M,g)$ is isometric to the Euclidean space $(\mathbb R^n,g_0),$ and $ \Omega_t\subset M$ is isometric to the ball $\Omega_t^\star\subset \mathbb R^n$  	
 	for a.e. $0<t<\|u\|_{L^\infty(M)}$. 
 \end{proposition}

 

 {\it Proof.} 
  The proof is inspired by Aubin \cite{Aubin}, and Brothers and Ziemer \cite{BZ}. Without loss of generality, we may assume that $u\geq 0$ 
 since $\|\nabla_g u\|_{L^p(M)}=\|\nabla_g |u|\|_{L^p(M)}$, and by density, it it enough to consider functions belonging to $\mathcal{C}^{n}(M).$ For every $0<t<\|u\|_{L^\infty(M)}=:L$, let
 $$\Pi_t:=u^{-1}(t)\subset M\ \ {\rm  and}\ \ \Pi_t^\star:=(u^\star)^{-1}(t)\subset \mathbb R^n,$$
 and due to \eqref{omegak}, 
  \begin{align*}\mathcal V(t):= {\rm Vol}_g(\Omega_t)={\rm Vol}_{g_0}(\Omega_t^\star).
 \end{align*}
  Consider the set of critical points $C=\{x\in M:\nabla_g u(x)=0\}$ and $C^\star=\{x\in \mathbb R^n:\nabla u^\star(x)=0\},$ respectively. Since $u\in \mathcal{C}^{n}(M)$, we have that the set of critical values $u(C)$ is a null measure set in $\mathbb{R}$ by Sard's theorem. The set $\Pi_t$ for $t \notin u(C)$ is a smooth, regular surface  of class $\mathcal{C}^{n}$. Similarly we have that $u^{\star}(C^{\star})$ is also of null measure.

 The co-area formula  implies that
 \begin{eqnarray}
 	\mathcal V(t)&=&{\rm Vol}_g(C\cap u^{-1}(t,L))+\int_t^L\left(\int_{\Pi_s}\frac{1}{|\nabla_g u|} {\rm d}\mathcal H^{n-1}\right){\rm d}s\label{first-repr}\\&=&{\rm Vol}_{g_0}(C^\star\cap (u^\star)^{-1}(t,L))+\int_t^L\left(\int_{\Pi_s^\star}\frac{1}{|\nabla u^\star|} {\rm d}\mathcal H^{n-1}\right){\rm d}s.\label{second-repr} 
 \end{eqnarray}
%
 
 By the monotonicity of the function  $\mathcal V$, we conclude that it is differentiable almost everywhere. Furthermore notice that the functions $t\mapsto {\rm Vol}_g(C\cap u^{-1}(t,L))$, $ t \mapsto {\rm Vol}_{g_0}(C^\star\cap (u^\star)^{-1}(t,L))$ have the same properties. Moreover, these latter functions have vanishing derivatives almost everywhere. Thus  
  \eqref{first-repr} and \eqref{second-repr}  imply that 
 \begin{equation}\label{V-deriv-1}
 	-\mathcal V'(t) =  \int_{\Pi_t}\frac{1}{|\nabla_g u|}{\rm d}\mathcal H^{n-1} = -\int_{\Pi_t^\star}\frac{1}{|\nabla u^\star|}{\rm d}\mathcal H^{n-1}\ \ {\rm for\ a.e.}\ 0<t<L.
 \end{equation}
In the sequel, we consider only those values of $t>0$ for which $\mathcal V'(t)$ is well-defined and the above formula holds.
 
 Since $u^\star$ is
 radially symmetric, the set $\Pi_t^\star$ is an $(n-1)$-dimensional
 sphere. Furthermore, $|\nabla u^\star|_t:=|\nabla u^\star|$ is constant on the $(n-1)$-dimensional  sphere
 $\Pi_t^\star$ and by (\ref{V-deriv-1}) it follows that
 \begin{equation}\label{V-DER-MASIK}
 	\mathcal V'(t)=-\frac{\mathcal H^{n-1}(\Pi_t^\star)}{|\nabla u^\star|_t}\ \ {\rm for\ a.e.}\ 0<t<L.
 \end{equation}
 By  (\ref{V-deriv-1}) and H\"older's inequality 
 we infer that
 \begin{eqnarray*}
 	{\mathcal H^{n-1}}(\Pi_t)&=&\int_{\Pi_t}{{\rm d} \mathcal H^{n-1}}=\int_{\Pi_t}\frac{1}{|\nabla_gu|^\frac{p-1}{p}}|\nabla_gu|^\frac{p-1}{p}{\rm d}\mathcal H^{n-1} \\&\leq&
 	\left(\int_{\Pi_t}\frac{1}{|\nabla_g u|}{\rm d}\mathcal H^{n-1}\right)^\frac{p-1}{p}\left(\int_{\Pi_t}{|\nabla_g u|^{p-1}}{{\rm d} \mathcal H^{n-1}}\right)^\frac{1}{p}
 	\\&\leq&\left(-\mathcal V'(t)\right)^\frac{p-1}{p}\left(\int_{\Pi_t}{|\nabla_g u|^{p-1}}{{\rm d} \mathcal H^{n-1}}\right)^\frac{1}{p}.
 \end{eqnarray*}
 Since  ${\rm Vol}_g(\Omega_t)={\rm Vol}_{g_0}(\Omega_t^\star)$ and $\Pi_t$ is a $\mathcal{C}^n$ smooth regular surface for a.e. $0<t<L$  (and for such surfaces the $(n-1)$-dimensional Hausdorff measure and the perimeter coincide), by the isoperimetric inequality \eqref{eq-isoperimetric-1} we have 
 \begin{equation}\label{isoperi-2}
 	\mathcal H^{n-1}(\Pi_t)=\mathcal P_g(\Pi_t)\geq {\sf AVR}_g^\frac{1}{n}\ \mathcal P_{g_0}(\Pi_t^\star)= {\sf AVR}_g^\frac{1}{n}\ \mathcal H^{n-1}(\Pi_t^\star)\ \ {\rm for\ a.e.}\ 0<t<L.
 \end{equation}
 Thus, by  (\ref{isoperi-2}) and relation
 (\ref{V-DER-MASIK}), the previous estimate implies  that
 \begin{eqnarray}\label{iso-egyenloseg}
 	\nonumber 
 	\int_{\Pi_t}{|\nabla_g u|^{p-1}}{\rm d}
 		\mathcal H^{n-1} &\geq& \left(\mathcal H^{n-1}(\Pi_t)\right)^p \left(-\mathcal V'(t)\right)^{1-p}
 	\\&\geq&  {\sf AVR}_g^\frac{p}{n}\left(\mathcal H^{n-1}(\Pi_t^\star)\right)^p \left(\frac{\mathcal H^{n-1}(\Pi_t^\star)}{|\nabla u^\star|_t}\right)^{1-p}\\
 	&=& {\sf AVR}_g^\frac{p}{n}\int_{\Pi_t^\star}{|\nabla u^\star|^{p-1}}{\rm d}
 		\mathcal H^{n-1}.\nonumber
 \end{eqnarray}
 By combining again the co-area formula with this estimate, it follows that
 \begin{eqnarray*}\label{Polya-Szego}
 	\nonumber \int_{M}{|\nabla_g u|^{p}}{{\rm d}
 	}v_g &=& \int_0^\infty\int_{\Pi_t}{|\nabla_g u|^{p-1}}{\rm d}
 \mathcal H^{n-1} {\rm d}t \\
 	&\geq&{\sf AVR}_g^\frac{p}{n}\int_0^\infty \int_{\Pi_t^\star}{|\nabla
 		u^\star|^{p-1}}{\rm d}
 	\mathcal H^{n-1}{\rm d}t={\sf AVR}_g^\frac{p}{n}\int_{\mathbb
 		R^n}{|\nabla u^\star|^{p}}{\rm d}x,
 \end{eqnarray*}
 which concludes the  proof of inequality \eqref{harmadik-tul}.
 
 If equality holds in \eqref{harmadik-tul} for some nonnegative $u\in \mathcal{C}^{n}(M)\setminus \{0\}$, by  \eqref{isoperi-2} we necessarily have  for  a.e. $ 0<t<L$ 
  that 
 $$	\mathcal P_g(\Pi_t)= {\sf AVR}_g^\frac{1}{n}\ \mathcal P_{g_0}(\Pi_t^\star)=n\omega_n^\frac{1}{n} \ {\sf AVR}_g^\frac{1}{n}{\rm Vol}_{g_0}(\Omega_t^\star)^\frac{n-1}{n}=n\omega_n^\frac{1}{n} \ {\sf AVR}_g^\frac{1}{n}{\rm Vol}_g(\Omega_t)^\frac{n-1}{n}.$$
 Since by Sard's theorem, for almost every $t$ the isoperimetric domains $\Omega_t$ have regular $\mathcal{C}^n$ boundaries,  Theorem \ref{tetel-egyenloseg} can be applied. Thus 
 ${\sf AVR}_g=1,$ i.e., $(M,g)$ is isometric to the Euclidean space $(\mathbb R^n,g_0),$ and $ \Omega_t\subset M$ is isometric to the ball $\Omega_t^\star\subset \mathbb R^n$  	
 for a.e. $0<t<L$. 
  \hfill$\square$

 \subsection{Gagliardo-Nirenberg interpolation inequality} \label{sect1.1}
 Sharp Gagliardo-Niren\-berg ine\-qua\-lities on $\mathbb R^n$ are known after Del Pino and Dolbeault \cite{DelPino-Dolb} and Cordero-Erausquin, Nazaret and Villani \cite{CE-N-Villani}. In the sequel we establish a sharp Gagliardo-Nirenberg inequality on Riemannian manifolds with ${\sf Ric}\geq 0$ whose particular form provides the statement of 
 Theorem \ref{main-0}.

 When  $p\in(1,n)$  and  $1<\alpha\leq \frac{n}{n-p}$, the
  \textit{Gagliardo-Nirenberg inequality} on $\mathbb R^n$ reads as 
 \begin{equation}\label{Villani-1}
 	\|u\|_{L^{\alpha p}(\mathbb R^n)}\leq \mathcal G_{\alpha,p,n}
 	\|\nabla u\|_{L^p(\mathbb R^n)}^{\theta}\|u\|_{L^{\alpha(p-1)+1}(\mathbb R^n)}^{1-\theta},\
 	\forall u\in \dot W^{1,p}(\mathbb R^n),
 \end{equation}
 where 
 \begin{equation}\label{theta-best}
 	\theta=\frac{p^\star(\alpha-1)}{\alpha p(p^\star-\alpha
 		p+\alpha-1)},
 \end{equation}
 and the best constant 
 \begin{equation}\label{GN-1-konstans}
 	\mathcal
 	G_{\alpha,p,n}:=\left(\frac{\alpha-1}{p'}\right)^\theta
 	\frac{\left(\frac{p'}{n}\right)^{\frac{\theta}{p}+\frac{\theta}{n}}\left(\frac{\alpha (p-1)+1}{\alpha
 			-1}-\frac{n}{p'}\right)^\frac{1}{\alpha p}
 		\left(\frac{\alpha (p-1)+1}{\alpha
 			-1}\right)^{\frac{\theta}{p}-\frac{1}{\alpha p}}}{\left(\omega_n {\sf B}\left(\frac{\alpha (p-1)+1}{\alpha
 			-1}-\frac{n}{p'},\frac{n}{p'}\right)\right)^{\frac{\theta}{n}}}
 \end{equation}
 is achieved by the unique family of functions
 $h_{\alpha,p}^\lambda(x)=(\lambda+|x|^\frac{p}{p-1})^\frac{1}{1-\alpha},\ x\in \mathbb
 R^n,$
 $\lambda>0.$ Here, $p^\star=\frac{np}{n-p}$, $p'=\frac{p}{p-1}$ and $\dot W^{1,p}(\mathbb R^n)=\{u\in L^{p^\star}(\mathbb R^n):|\nabla u|\in L^p(\mathbb
 R^n)\},$ while   ${\sf B}(\cdot,\cdot)$ is the Euler
 beta-function.
 Let us recall from the Introduction the function space $$\dot W^{1,p}(M)=\{u\in L^{p^\star}(M):|\nabla_g u|\in L^p(M)\}.$$ 
 
 
 \begin{theorem}\label{main-GN} Let $(M,g)$ be an $n$-dimensional Riemannian manifold as in Theorem \ref{main-0}.   Let $p\in (1,n)$, $1<\alpha\leq \frac{n}{n-p}$ and the constants $\theta$ and $\mathcal G_{\alpha,p,n}$ given by \eqref{theta-best} and \eqref{GN-1-konstans}, respectively.   Then  
 	\begin{equation}\label{egyenlet-1}
 		\|u\|_{L^{\alpha p}(M)}\leq {\sf K}_g^{\sf GN}\
 		\|\nabla_g u\|_{L^p(M)}^{\theta}\|u\|_{L^{\alpha(p-1)+1}(M)}^{1-\theta},\ \ \forall u\in \dot W^{1,p}(M),
 	\end{equation}
 	where the constant $ {\sf K}_g^{\sf GN}=\mathcal G_{\alpha,p,n}\ {\sf AVR}_g^{-\frac{\theta}{n}}$ is  sharp. Moreover,  equality holds in \eqref{egyenlet-1} for some nonzero and nonnegative function $u\in \mathcal{C}^{n}(M)\cap 
 	\dot W^{1,p}(M)$
 	if and only if ${\sf AVR}_g=1$ and $u=h_{\alpha,p}^\lambda$ a.e. for some $\lambda>0$ $($up to isometry$).$
 \end{theorem}
 
 {\it Proof.} Clearly, it is enough to prove \eqref{egyenlet-1} for nonnegative functions $u\in \mathcal{C}^{\infty}_0(M)$; the inequality for general $u \in \dot W^{1,p}(M)$ will follow by approximation.  
 We recall that the Euclidean rearrangement function $u^\star$ of $u$
 satisfies the optimal Gagliardo-Nirenberg inequality
 (\ref{Villani-1}), thus relations \eqref{masodik-tul} and \eqref{harmadik-tul} imply that
 \begin{eqnarray}\label{g-egyenloseg}
 	\nonumber 
 	\|u\|_{L^{\alpha p}(M)} &= & \|u^\star\|_{L^{\alpha p}(\mathbb R^n)} 
 	\\&	\leq&  \mathcal G_{\alpha,p,n}
 	\|\nabla u^\star\|_{L^p(\mathbb R^n)}^{\theta}\|u^\star\|_{L^{\alpha(p-1)+1}(\mathbb R^n)}^{1-\theta} \\
 	&\leq&\mathcal G_{\alpha,p,n}\ {\sf AVR}_g^{-\frac{\theta}{n}}  \|\nabla_g
 	u\|_{L^p(M)}^{\theta}\|u\|_{L^{\alpha(p-1)+1}(M)}^{1-\theta},\nonumber
 \end{eqnarray}
 which proves inequality \eqref{egyenlet-1}. 
 
 To prove the statement about sharpness, assume by contradiction, that the constant $ {\sf K}_g^{\sf GN}=\mathcal G_{\alpha,p,n}\ {\sf AVR}_g^{-\frac{\theta}{n}}$ is not sharp in \eqref{egyenlet-1}, i.e., there exists $0< C<{\sf K}_g^{\sf GN}$  such that 
 $$\|u\|_{L^{\alpha p}(M)}\leq  C\
 \|\nabla_g u\|_{L^p(M)}^{\theta}\|u\|_{L^{\alpha(p-1)+1}(M)}^{1-\theta},\
 \forall u\in \dot W^{1,p}(M).$$
 Since $(M,g)$ is a complete Riemannian manifold with ${\sf Ric}\geq 0$, the validity of the latter inequality implies the quantitative non-collapsing volume property
 $${\rm Vol}_g(B_{x}(r))\geq\left(\frac{\mathcal G_{\alpha,p,n}}{ C}\right)^\frac{n}{\theta}\omega_n r^n\ \  {\rm for\ all}\ x\in M\ {\rm and}\
 r\geq 0,$$ 
 see Ledoux \cite{Ledoux-CAG} (for $\theta=1$) and  Xia \cite{Xia-JFA}, Krist\'aly \cite{Kristaly-Calculus}, and Krist\'aly and Ohta \cite{Kri-Ohta} (for general $\theta$ from \eqref{theta-best}). 
 By the latter relation, one has that
 $${\sf AVR}_g=\lim_{r\to \infty}\frac{{\rm Vol}_g(B_x(r))}{\omega_n r^n}\geq \left(\frac{\mathcal G_{\alpha,p,n}}{ C}\right)^\frac{n}{\theta}.$$
 This inequality is equivalent to $ C\geq {\sf K}_g^{\sf GN}$ which  contradicts our initial assumption $0< C<{\sf K}_g^{\sf GN}$. In conclusion, the constant ${\sf K}_g^{\sf GN}$ is sharp in inequality \eqref{egyenlet-1}.
 
 Assume now that equality holds in \eqref{egyenlet-1} for some nonzero and nonnegative function $u\in \mathcal{C}^{n}(M)\cap \dot W^{1,p}(M)$. Clearly, $u$ is fast decaying and  by \eqref{masodik-tul}, \eqref{harmadik-tul} and \eqref{Villani-1} we have
 \begin{eqnarray*}
 	\nonumber 
 	\|u\|_{L^{\alpha p}(M)} &= & {\sf K}_g^{\sf GN}  \|\nabla_g
 	u\|_{L^p(M)}^{\theta}\|u\|_{L^{\alpha(p-1)+1}(M)}^{1-\theta}\geq
 	{\sf K}_g^{\sf GN}{\sf AVR}_g^{\frac{\theta}{n}}  \|\nabla
 	u^\star\|_{L^p(\mathbb R^n)}^{\theta}\|u^\star\|_{L^{\alpha(p-1)+1}(\mathbb R^n)}^{1-\theta}\\&\geq&
 	\|u^\star\|_{L^{\alpha p}(\mathbb R^n)} 
 	=  \|u\|_{L^{\alpha p}(M)}.
 \end{eqnarray*}
 Consequently, we have equalities in the above chain of inequalities. In particular, by the first equality we have equality in the P\'olya-Szeg\H o inequality \eqref{harmadik-tul}; thus Proposition \ref{PSz-proposition} implies that ${\sf AVR}_g=1$, i.e., $(M,g)$ is isometric to the Euclidean space $(\mathbb R^n,g_0).$ Furthermore, the second equality implies that $u^\star$ is an extremal in the Euclidean Gagliardo-Nirenberg inequality \eqref{Villani-1}, which shows that $u^\star=h_{\alpha,p}^\lambda$ for some $\lambda>0$. Since the set $$\{x\in \mathbb R^n:u^\star(x)>0,\ \nabla u^\star(x)=0\}$$ is a singleton containing only $0\in \mathbb R^n$, it turns out by  Brothers and Ziemer \cite[Theorem 1.1]{BZ} that  (up to an isometry), 
  $u=u^\star$ a.e. 
 The converse is trivial. \hfill $\square$


%

\begin{remark}\rm 
	(i) Theorem \ref{main-0} is a simple consequence of Theorem \ref{main-GN} by considering $\theta=1.$
	
	(ii) Further sharp Sobolev-type inequalities can be proved on Riemannian manifolds with 
	with ${\sf Ric}\geq 0$, similarly to Theorem  
	\ref{main-GN}, as the dual Gagliardo-Nirenberg inequality ($0<\alpha<1$), $L^p$-log-Sobolev inequality $(\alpha \to 1)$, and Faber-Krahn inequality $(\alpha \to 0)$, see \cite{CE-N-Villani}. In fact, these inequalities  follow by Proposition \ref{PSz-proposition}, while their sharpness by Krist\'aly \cite{Kristaly-Calculus}.
\end{remark}

\subsection{Rayleigh-Faber-Krahn inequality:  first eigenvalues in sharp form}\label{Poincare-section}

In this subsection we  prove Theorem \ref{RFK-1-0};  the following auxiliary result will be used whose proof is based on a simple application of the layer cake representation.  

\begin{lemma}\label{lemma-layer}
	Let $(M,g)$ be a complete Riemannian manifold, $R>0$ and $x_0\in M$ be arbitrarily fixed,  and $f:[0,R]\to \mathbb R$ be a $\mathcal C^1$-function on $(0,R)$. Then 
	$$	\int_{B_{x_0}(R)}f(d_g(x_0,x)){\rm d}v_g=f(R){\rm Vol}_g(B_{x_0}(R))-\int_0^R f'(r){\rm Vol}_g(B_{x_0}(r)){\rm d}r.
	$$
\end{lemma} 


The precursor of Theorem \ref{RFK-1-0} reads as follows.


\begin{theorem}\label{main-1} Let $(M,g)$ be an $n$-dimensional Riemannian manifold as in Theorem \ref{main-0}.   Then for every smooth bounded open set $\Omega\subset M$ and $ u\in W_0^{1,2}(\Omega)$   we have 
	\begin{equation}\label{egyenlet-poincare}
		\displaystyle  {{{\rm Vol}_g(\Omega)^{-\frac{2}{n}}} \int_{\Omega} u^2{\rm d}v_g\leq
			{{\sf R}_g}\, 
			\int_{\Omega} |\nabla_g u|^2 {\rm d}v_g},
	\end{equation}
	where the constant $ {\sf R}_g= j^{-2}_{\frac{n}{2}-1}(\omega_n{\sf AVR}_g)^{-\frac{2}{n}}$ is  sharp.  In addition, equality holds in \eqref{egyenlet-poincare} for some set $\Omega\subset M$ and for some nonzero  function $u\in W_0^{1,2}(M)$ if and only if $(M,g)$ is  isometric to $(\mathbb R^n,g_0)$ and $\Omega$ is isometric to a ball $B\subset \mathbb R^n$.
\end{theorem}


{\it Proof.}  Let $\Omega\subset M$ be any smooth bounded open set and $B\subset \mathbb R^n$ be a ball with ${\rm Vol}_g(\Omega)={\rm Vol}_{g_0}(B)$.  If $u:\Omega\to \mathbb R$ is any nonzero function with the usual smoothness properties, its 
Euclidean rearrangement function $u^\star:B\to
[0,\infty)$ satisfies the properties \eqref{elso-tul}, \eqref{masodik-tul} and \eqref{harmadik-tul}; in particular, the latter relations combined with the Euclidean Rayleigh-Faber-Krahn inequality immediately gives 
\begin{equation}\label{BVP}
\frac{\ds\int_{\Omega} |\nabla_g u|^2 {\rm d}v_g}{{{\rm Vol}_g(\Omega)^{-\frac{2}{n}}} \ds\int_{\Omega} u^2{\rm d}v_g}\geq \frac{ {\sf AVR}_g^\frac{2}{n} \ds\int_{B} |\nabla u^\star|^2 {\rm d}x}{{{\rm Vol}_{g_0}(B)^{-\frac{2}{n}}} \ds\int_{B} (u^\star)^2{\rm d}x}\geq j^{2}_{\frac{n}{2}-1}(\omega_n{\sf AVR}_g)^{\frac{2}{n}}={\sf R}_g^{-1},
\end{equation}
which is exactly  inequality \eqref{egyenlet-poincare}.

The sharpness of the constant ${\sf R}_g= j^{-2}_{\frac{n}{2}-1}({\omega_n\sf AVR}(g))^{-\frac{2}{n}}$ in \eqref{egyenlet-poincare} is more delicate, which requires fine properties of Bessel functions of the first kind $J_\nu$, $\nu\in \mathbb R,$ see e.g.\ Olver {\it et al.} \cite{Digital}; for completeness, we outline its proof. 
%
%
%
By contraction, we assume there exists $ C<{\sf R}_g$ such that for every smooth bounded open set $\Omega\subset M$, 
	\begin{equation}\label{egyenlet-poincare-2}
	\displaystyle  {\rm Vol}_g(\Omega)^{-\frac{2}{n}} \int_{\Omega} u^2{\rm d}v_g\leq  C \, \int_{\Omega} |\nabla_g u|^2 {\rm d}v_g,\ \ \forall u\in W_0^{1,2}(\Omega).
\end{equation}
For convenience of notation, we choose $\nu=\frac{n}{2}-1\geq 0.$ For every $R>0$ and $x_0\in M$, we consider the function $u_R:B_{x_0}(R)\to \mathbb R$ defined by 
 $$u_R(x)=d_g(x_0,x)^{-\nu}J_\nu\left(j_\nu\frac{d_g(x_0,x)}{R}\right),\ \ x\in B_{x_0}(R).$$
 It is clear  that \eqref{egyenlet-poincare-2} can be applied to the function $u_R$ and to the set $\Omega=B_{x_0}(R)$, i.e., 
 \begin{equation}\label{egyenlet-poincare-3}
 	\displaystyle  {\rm Vol}_g(B_{x_0}(R))^{-\frac{2}{n}} \int_{B_{x_0}(R)} u_R^2{\rm d}v_g \leq \mathcal C \, \int_{B_{x_0}(R)} |\nabla_g u_R|^2 {\rm d}v_g.
 \end{equation}
Basic properties of Bessel functions combined with the eikoinal equation 
$$|\nabla_g d_g(x_0,x)|=1\ \ {\rm for\ a.e.}\ x\in M,$$  imply that 
%
%
$$
	\displaystyle  \int_{B_{x_0}(R)} |\nabla_g u_R|^2 {\rm d}v_g= \frac{j_\nu^2}{R^2}\int_{B_{x_0}(R)}d_g(x_0,x)^{-2\nu} J_{\nu+1}^2\left(j_\nu\frac{d_g(x_0,x)}{R}\right) {\rm d}v_g.
$$
Applying Lemma \ref{lemma-layer}, a change of variables, the Lebesgue's dominated convergence theorem and finally an integration by parts, one has
\begin{eqnarray*}
	\lim_{R\to \infty}\int_{B_{x_0}(R)} |\nabla_g u_R|^2 {\rm d}v_g&=&j_\nu^2\omega_n{\sf AVR}_g\left(J_{\nu+1}^2(j_\nu)-\int_0^1t^n\frac{\rm d}{{\rm d}t}(t^{-2\nu}J_{\nu+1}^2(j_\nu t)){\rm d}t\right)\\&=&j_\nu^2n\omega_n{\sf AVR}_g\int_0^1tJ_{\nu+1}^2(j_\nu t)){\rm d}t.
\end{eqnarray*}
Since $J_\nu(j_\nu)=0$, a similar reasoning as above shows that 
\begin{eqnarray*}
	\lim_{R\to \infty}\frac{\ds\int_{B_{x_0}(R)}  u_R^2 {\rm d}v_g}{R^2}=n\omega_n{\sf AVR}_g\int_0^1tJ_{\nu}^2(j_\nu t)){\rm d}t.
\end{eqnarray*}
Letting $R\to \infty$ in \eqref{egyenlet-poincare-3} and taking into account the latter relations, we obtain that
$$(\omega_n{\sf AVR}_g)^{-\frac{2}{n}} \leq  C j_\nu^2.$$
This inequality clearly contradicts $ C<{\sf R}_g= j^{-2}_{\nu}({\omega_n\, \sf AVR}(g))^{-\frac{2}{n}},$  concluding the sharpness of  ${\sf R}_g$ in \eqref{egyenlet-poincare}.

Assume now that equality holds in \eqref{egyenlet-poincare} for some open bounded set $\Omega\subset M$ and some nonzero and nonnegative function $u\in W_0^{1,2}(\Omega)$. 
In particular, equality in \eqref{egyenlet-poincare} (thus in   \eqref{BVP}) implies equality in \eqref{harmadik-tul}, where $u^\star:B\to \mathbb R$ is the Euclidean rearrangement function of $u$, and $B\subset \mathbb R^n$ is a ball with   ${\rm Vol}_g(\Omega)={\rm Vol}_{g_0}(B)$. Furthermore, since $u$ verifies a second order PDE involving the   operator $\Delta_g$ (obtained as the Euler-Lagrange equation associated with \eqref{faber-krahn-0}), it turns out by standard regularity theory that $u\in \mathcal C^\infty(\Omega)$.  Consequently, by Proposition \ref{PSz-proposition} one has ${\sf AVR}_g=1$ and $\Omega$ is isometric to  $B\subset \mathbb R^n$.  The converse is trivial. 
\hfill
$\square$
\vspace{0.2cm}

{\it Proof of Theorem \ref{RFK-1-0}}. The proof rests upon Theorem \ref{main-1}. Indeed, the sharp estimate of the first eigenvalue  $\lambda_{1,g}^D(\Omega)$ in \eqref{egyenlet-poincare-1-0} immediately follows by \eqref{egyenlet-poincare}, noticing that ${\sf R}_g={\sf \Lambda}_g^{-1}.$ Standard compactness and variational argument show that the infimum in \eqref{faber-krahn-0} is achieved. Thus, if equality holds in \eqref{egyenlet-poincare-1-0} for some $\Omega\subset M$, then there exists  some nonzero and sign-preserving function $u\in W_0^{1,2}(\Omega)$ which produces equality in \eqref{egyenlet-poincare} as well, concluding the proof, cf. Theorem \ref{main-1}.  \hfill
$\square$

\begin{remark}\rm 
	Having sharp Sobolev inequalities on Riemannian manifolds with ${\sf Ric}\geq 0$ (see  \S \ref{sect1.1}-\ref{Poincare-section}), a natural question would be to establish such inequalities in the framework of general ${\sf CD}(0,N)$ spaces. Based on the sharp isoperimetric inequality in Theorem  \ref{main-isop-CD}, a P\'olya-Szeg\H o inequality is expected to hold on ${\sf CD}(0,N)$ spaces; a similar argument has been completed recently on ${\sf CD}(\kappa,N)$ 
	spaces with $\kappa>0$, see Mondino and Semola \cite{MSemola}.
\end{remark}

\noindent {\bf Acknowledgments.} The authors thank S. Brendle, G. Carron,  L. Mazzieri, E. Milman and S. Ohta   for stimulating conversations in the early stage of the manuscript. We also would like to thank the Reviewers for their thoughtful comments and efforts towards improving the presentation of our paper. 
 \\

\noindent \textbf{Final statement.} The authors state that there is no conflict of interest.


\begin{thebibliography}{99}
	
	


\bibitem{AFM} V. Agostiniani, M. Fogagnolo, L. Mazzieri, Sharp geometric inequalities for closed hypersurfaces in manifolds with nonnegative Ricci curvature. \textit{Invent. Math.} 222 (2020), no. 3, 1033--1101.

\bibitem{AFP}  L. Ambrosio, N. Fusco, D. Pallara,  Functions of bounded variation and free discontinuity problems. Oxford Mathematical Monographs. The Clarendon Press, Oxford University Press, New York, 2000.

\bibitem{Anderson} M. Anderson, On the topology of complete manifold of nonnegative
Ricci curvature. \textit{Topolog}y 3  (1990), 41--55.

\bibitem{Aubin} T. Aubin, Probl\`emes isop\'erim\'etriques et espaces de Sobolev.  \textit{J. Differential Geometry} 11 (1976), no. 4, 573--598.


\bibitem{BS} K. Bacher, K.-T. Sturm, Ricci bounds for Euclidean and spherical cones. \textit{Singular phenomena and scaling in mathematical models, 3-23, Springer, Cham.} (2014)

	\bibitem{BK-Adv} Z. M. Balogh, A. Krist\'{a}ly, Equality in Borell-Brascamp-Lieb inequalities on curved spaces. \textit{Adv. Math.} 339 (2018) 453-494





\bibitem{Brendle} S. Brendle, Sobolev inequalities in manifolds with nonnegative curvature. \textit{Comm. Pure. Appl. Math.}, to appear. 


\bibitem{BZ}  J. E. Brothers, W. P. Ziemer,  Minimal rearrangements of Sobolev functions. \textit{J. Reine Angew. Math.} 384 (1988), 153--179.



\bibitem{Cabre-Ros-Oton}  X. Cabr\'e, X. Ros-Oton,  Sobolev and isoperimetric inequalities with monomial weights. \textit{J. Differential Equations} 255 (2013), no. 11, 4312--4336. 

\bibitem{Cabre-Ros-Oton-Serra}  X. Cabr\'e, X. Ros-Oton, J. Serra, Sharp isoperimetric inequalities via the ABP method. \textit{J. Eur. Math. Soc.} 18 (2016),  2971--2998. 

 
\bibitem{Carron} G. Carron, Euclidean volume growth for complete Riemannian manifolds.  \textit{Milan J. Math.} 88 (2020), 455--478.

\bibitem{Carron2} G. Carron, Inegalit\'es isop\'erim\'etriques de Faber-Krahn et consequences. Actes de la Table Ronde de G\'eom\'etrie Diff\'erentielle (Luminy, 1992), 205--232, S\'emin. Congr., 1, Soc. Math. France, Paris, 1996. 

\bibitem{Cavalletti-Mondino} F. Cavalletti, A. Mondino,
Sharp and rigid isoperimetric inequalities in metric-measure spaces with lower Ricci curvature bounds. 
\textit{Invent. Math.} 208 (2017), no. 3, 803--849. 


\bibitem{Chavel} I. Chavel,  Riemannian Geometry. A Modern Introduction.
Second Edition, 2006.


\bibitem{CC} J. Cheeger, T. H.  Colding, 
Lower bounds on Ricci curvature and the almost rigidity of warped products.
\textit{Ann. of Math. (2)}  144 (1996), no. 1, 189--237. 








\bibitem{Cinti-etal} E. Cinti, F. Glaudo, A. Pratelli, X. Ros-Oton, J. Serra, Sharp quantitative stability for isoperimetric inequalities with homogeneous weights. \textit{Trans. Amer. Math. Soc.}, accepted (2021). 


\bibitem{Colding} T. H. Colding, 
Ricci curvature and volume convergence.
\textit{Ann. of Math. (2)}  145 (1997), no. 3, 477--501. 

	\bibitem{CEMS} D. Cordero-Erausquin, R. J. McCann, M. Schmuckenschl\"{a}ger, 
A Riemannian interpolation inequality \`{a} la Borell, Brascamp and Lieb. \textit{Invent. Math.} 146, (2001), 219--257.

\bibitem{CE-N-Villani}  D. Cordero-Erausquin, B. Nazaret, C. Villani, A mass-transportation approach to sharp Sobolev and Gagliardo-Nirenberg inequalities. \textit{Adv. Math.} 182 (2004), no. 2, 307--332.




\bibitem{Coulhon-SC}  T. Coulhon, L. Saloff-Coste,  Isop\'erim\'etrie pour les groupes et les vari\'et\'es.  \textit{Rev. Mat. Iberoamericana} 9 (1993), no. 2, 293--314.

\bibitem{Croke} C. Croke, A sharp four-dimensional isoperimetric inequality.
\textit{Comment. Math. Helv.} 59 (1984), no. 2, 187--192.

\bibitem{DelPino-Dolb} M. Del Pino, J. Dolbeault, The optimal Euclidean $L^p$-Sobolev logarithmic inequality. \textit{J. Funct. Anal.} 197 (2003)
151--161.




\bibitem{Druetetal} O. Druet, E. Hebey, M. Vaugon, Optimal Nash's inequalities on
Riemannian manifolds: the influence of geometry.  \textit{Int. Math. Res. Not.} IMRN, 14, 1999, 735--779.



\bibitem{doCarmo-Xia} M. P. do Carmo, C. Xia,
Complete manifolds with non-negative Ricci curvature and the
Caffarelli-Kohn-Nirenberg inequalities. \textit{Compos. Math.} {140} (2004),
818--826.

  

\bibitem{Fogagnolo-Mazzieri} M. Fogagnolo, L. Mazzieri, 
Minimising hulls, $p$-capacity and isoperimetric inequality on complete Riemannian manifolds.  Preprint, December 2020. Link: https://arxiv.org/abs/2012.09490

\bibitem{Fusco} N. Fusco, The classical isoperimetric theorem. Link: https://www.docenti.unina.it/webdocenti-be/allegati/materiale-didattico/377226


\bibitem{Gentil} I. Gentil, The general optimal $L^p$-Euclidean logarithmic Sobolev inequality by Hamilton-Jacobi equations. \textit{J. Funct. Anal.} 202 (2003), no. 2, 591--599.

\bibitem{GW1} R. E. Greene, H. Wu, $\mathcal{C}^{\infty} $ Convex functions and manifolds of positive curvature. \textit{ Acta Math.} 137 (1976),  209--245.
\bibitem{GW} R. E. Greene, H. Wu, $\mathcal{C}^{\infty} $ approximations of convex, subharmonic, and plurisubharmonic functions. \textit{Ann. Sci. ENS} 12, 1 (1979), 48--84.




\bibitem{Hebey} E. Hebey, Nonlinear analysis on manifolds: Sobolev spaces and inequalities. Courant Lecture Notes in Mathematics, 5.
New York University, Courant Institute of Mathematical Sciences, New
York; American Mathematical Society, Providence, RI, 1999.

\bibitem{Johne} F. Johne, Sobolev inequalities on manifolds with nonnegative Bakry-\'Emery Ricci curvature, preprint, March 2021. Link: https://arxiv.org/abs/2103.08496.

\bibitem{Ketterer} C. Ketterer,  Cones over metric measure spaces and the maximal diameter theorem.\textit{ J. Math. Pures Appl.}  (9) 103 (2015), no. 5, 1228–1275.

\bibitem{Klartag} B. Klartag, Needle decompositions in Riemannian geometry. \textit{ Mem. Amer. Math. Soc.}
 249 (2017), no. 1180, 

%
\bibitem{Kleiner} B. Kleiner, An isoperimetric comparison theorem. \textit{Invent. Math.} 108 (1992), no. 1, 37--47.
%



\bibitem{Kolesnikov-Milman} A.V. Kolesnikov, E. Milman, Poincar\'e and Brunn-Minkowski inequalities on the boundary of weighted Riemannian manifolds. \textit{Amer. J. Math. } 140 (2018), no. 5, 1147--1185.

\bibitem{Kristaly-Calculus} A. Krist\'aly,  Metric measure spaces supporting Gagliardo-Nirenberg inequalities: volume non-collapsing and rigidities. \textit{Calc. Var. Partial Differential Equations} 55 (2016), no. 5, Art. 112, 27 pp.

 
 
 
 \bibitem{Kri-Ohta} A. Krist\'aly, S. Ohta, Caffarelli-Kohn-Nirenberg inequality on metric measure spaces with applications. \textit{Math. Ann.} 357 (2013), no. 2, 711--726.
 





\bibitem{Ledoux-CAG} M. Ledoux, On manifolds with non-negative Ricci curvature and Sobolev inequalities.
\textit{Comm. Anal. Geom.} 7 (1999), no. 2, 347--353.

\bibitem{LL} E. H. Lieb, M. Loss, Analysis. Second edition. Graduate Studies in Mathematics, 14. American Mathematical Society, Providence, RI, 2001.

\bibitem{Liu} G. Liu, 3-manifolds  with  nonnegative  Ricci  curvature. \textit{Invent. Math.} 193 (2013), no. 2, 367--375.


\bibitem{Li} P. Li, Large time behavior of the heat equation on complete manifolds with nonnegative Ricci curvature. \textit{Ann. of Math. (2)}  124 (1986), no. 1, 1--21.

\bibitem{LV}  J. Lott, C. Villani, Ricci curvature for metric measure spaces via optimal transport.\textit{ Ann. of Math. (2)}  169 (3)
(2009),
903--991.
%


\bibitem{McCann} R. J. McCann, Polar factorization of maps on Riemannian manifolds. \textit{Geom. Funct. Anal.} 11 (3) (2001) 589--608. 

\bibitem{Menguy} X. Menguy,  Noncollapsing examples with positive Ricci curvature and infinite topological type. \textit{Geom. Funct. Anal.} 10 (2000), no. 3, 600--627.

\bibitem{Milman} E. Milman, Sharp isoperimetric inequalities and model spaces for the curvature-dimension-diameter condition. \textit{J. Eur. Math. Soc}. 17 (5) (2015) 1041--1078. 

\bibitem{Milman2} E. Milman,  Beyond traditional curvature-dimension I: new model spaces for isoperimetric and concentration inequalities in negative dimension. \textit{Trans. Amer. Math. Soc.} 369 (2017), no. 5, 3605--3637.

\bibitem{Milman-Rotem} E. Milman, L.  Rotem, Complemented Brunn-Minkowski inequalities and isoperimetry for homogeneous and non-homogeneous measures. \textit{Adv. Math. } 262 (2014), 867--908.

\bibitem{MSemola} A. Mondino, D. Semola, 
Polya-Szego inequality and Dirichlet $p$-spectral gap for non-smooth spaces with Ricci curvature bounded below. 
\textit{J. Math. Pures Appl.} (9) 137 (2020), 238--274. 


\bibitem{MS} M. Muratori, N. Soave, Some rigidity results for Sobolev inequalities and related PDEs on Cartan-Hadamard manifolds. Preprint: March 2021. Link: https://arxiv.org/abs/2103.08240
 


%


\bibitem{Ohta-1}  S. Ohta, Needle decompositions and isoperimetric inequalities in Finsler geometry.
\textit{J. Math. Soc. Japan} 70 (2018), 651--693.

\bibitem{Ohta-2}  S. Ohta, A semigroup approach to Finsler geometry: Bakry-Ledoux's isoperimetric inequality.  \textit{Comm. Anal. Geom.}, to appear. 


%
%

\bibitem{Digital} F. W. J. Olver, D. W. Lozier, R. F. Boisvert, C. W. Clark (eds.), NIST Handbook
of Mathematical Functions, Cambridge University Press, Cambridge, 2010.

\bibitem{Perelman-JAMS} G. Perelman, Manifolds of positive Ricci curvature with almost maximal volume. \textit{J. Am. Math. Soc.} 7 (1994),
299--305.

\bibitem{Petersen} P. Petersen, Riemannian geometry, Third edition, Graduate Texts in Mathematics, Springer, 2016.

\bibitem{Reiris} M. Reiris,  On Ricci curvature and volume growth in dimension three. \textit{J. Differential Geom.} 99 (2015), no. 2, 313--357.



%
%
%


	\bibitem{Sturm-1} K.-T. Sturm, On the geometry of metric measure spaces. I, \textit{Acta
	Math.} 196 (1) (2006), 65--131.

\bibitem{Sturm-2} K.-T. Sturm, On the geometry of metric measure spaces. II, \textit{Acta
	Math.} 196 (1) (2006), 133--177.

\bibitem{Talenti} G. Talenti,  Best constant in Sobolev inequality.\textit{ Ann. Mat. Pura Appl.} 110 (1976),
353--372.

\bibitem{Talenti-2} G. Talenti, Inequalities in rearrangement invariant function spaces. Nonlinear Analysis, Function Spaces and Applications, vol. 5 (Prague, 1994), 177--230. Prometheus, Prague (1994). 

\bibitem{Villani}  C. Villani, Optimal transport. Old and new. Grundlehren der Mathematischen Wissenschaften, vol. 338. Springer-Verlag, Berlin, 2009.


	\bibitem{WZ} Y. Wang, X. Zhang, An Alexandroff-Bakelman-Pucci estimate on Riemannian manifolds. \textit{Adv. Math.} 232 (2013), 499--512.




\bibitem{Weil} A. Weil, Sur les surfaces a courbure n\'egative. \textit{C. R. Acad. Sci. Paris S\'er I Math}. 182
(1926), 1069--1071.

\bibitem{Wulff} G. Wulff,     Zur   Frage   der   Geschwindigkeit   des   Wachstums   und   der   Aufl\"osung   der Kristallfl\"achen. \textit{Zeitschrift f\"ur Kristallographie} 34 (1901), 449--530.



\bibitem{Xia-JFA} C. Xia, The Gagliardo-Nirenberg inequalities and manifolds of non-negative Ricci curvature. \textit{J. Funct. Anal.} 224 (2005), no. 1, 230--241.

\bibitem{Zhu} S.-H. Zhu, A finiteness theorem for Ricci curvature in dimension three. \textit{J. Differential Geom.} 37 (1993), no. 3, 711--727.







\end{thebibliography}
\end{document}